\documentclass[11pt]{article}
\long\def\ig#1{\relax}
\ig{Thanks to Roberto Minio for this def'n.  Compare the def'n of
\comment in AMSTeX.}

\newcount \coefa
\newcount \coefb
\newcount \coefc
\newcount\tempcounta
\newcount\tempcountb
\newcount\tempcountc
\newcount\tempcountd
\newcount\xext
\newcount\yext
\newcount\xoff
\newcount\yoff
\newcount\gap%
\newcount\arrowtypea
\newcount\arrowtypeb
\newcount\arrowtypec
\newcount\arrowtyped
\newcount\arrowtypee
\newcount\height
\newcount\width
\newcount\xpos
\newcount\ypos
\newcount\run
\newcount\rise
\newcount\arrowlength
\newcount\halflength
\newcount\arrowtype
\newdimen\tempdimen
\newdimen\xlen
\newdimen\ylen
\newsavebox{\tempboxa}%
\newsavebox{\tempboxb}%
\newsavebox{\tempboxc}%

\makeatletter
\def\settypes(#1,#2,#3){\arrowtypea#1 \arrowtypeb#2 \arrowtypec#3}
\def\settoheight#1#2{\setbox\@tempboxa\hbox{#2}#1\ht\@tempboxa\relax}%
\def\settodepth#1#2{\setbox\@tempboxa\hbox{#2}#1\dp\@tempboxa\relax}%
\def\settokens[#1`#2`#3`#4]{%
     \def\tokena{#1}\def\tokenb{#2}\def\tokenc{#3}\def\tokend{#4}}
\def\setsqparms[#1`#2`#3`#4;#5`#6]{%
\arrowtypea #1
\arrowtypeb #2
\arrowtypec #3
\arrowtyped #4
\width #5
\height #6
}
\def\setpos(#1,#2){\xpos=#1 \ypos#2}

\def\bfig{\begin{picture}(\xext,\yext)(\xoff,\yoff)}
\def\efig{\end{picture}}

\def\putbox(#1,#2)#3{\put(#1,#2){\makebox(0,0){$#3$}}}

\def\settriparms[#1`#2`#3;#4]{\settripairparms[#1`#2`#3`1`1;#4]}%

\def\settripairparms[#1`#2`#3`#4`#5;#6]{%
\arrowtypea #1
\arrowtypeb #2
\arrowtypec #3
\arrowtyped #4
\arrowtypee #5
\width #6
\height #6
}

\def\resetparms{\settripairparms[1`1`1`1`1;500]\width 500}%default values%

\resetparms

\def\mvector(#1,#2)#3{%%
\put(0,0){\vector(#1,#2){#3}}%
\put(0,0){\vector(#1,#2){30}}%
}
\def\evector(#1,#2)#3{{%%
\arrowlength #3
\put(0,0){\vector(#1,#2){\arrowlength}}%
\advance \arrowlength by-30
\put(0,0){\vector(#1,#2){\arrowlength}}%
}}

\def\horsize#1#2{%
\settowidth{\tempdimen}{$#2$}%
#1=\tempdimen
\divide #1 by\unitlength
}

\def\vertsize#1#2{%
\settoheight{\tempdimen}{$#2$}%
#1=\tempdimen
\settodepth{\tempdimen}{$#2$}%
\advance #1 by\tempdimen
\divide #1 by\unitlength
}

\def\vertadjust[#1`#2`#3]{%
\vertsize{\tempcounta}{#1}%
\vertsize{\tempcountb}{#2}%
\ifnum \tempcounta<\tempcountb \tempcounta=\tempcountb \fi
\divide\tempcounta by2
\vertsize{\tempcountb}{#3}%
\ifnum \tempcountb>0 \advance \tempcountb by20 \fi
\ifnum \tempcounta<\tempcountb \tempcounta=\tempcountb \fi
}

\def\horadjust[#1`#2`#3]{%
\horsize{\tempcounta}{#1}%
\horsize{\tempcountb}{#2}%
\ifnum \tempcounta<\tempcountb \tempcounta=\tempcountb \fi
\divide\tempcounta by20
\horsize{\tempcountb}{#3}%
\ifnum \tempcountb>0 \advance \tempcountb by60 \fi
\ifnum \tempcounta<\tempcountb \tempcounta=\tempcountb \fi
}

\ig{ In this procedure, #1 is the paramater that sticks out all the way,
#2 sticks out the least and #3 is a label sticking out half way.  #4 is
the amount of the offset.}

\def\sladjust[#1`#2`#3]#4{%
\tempcountc=#4
\horsize{\tempcounta}{#1}%
\divide \tempcounta by2
\horsize{\tempcountb}{#2}%
\divide \tempcountb by2
\advance \tempcountb by-\tempcountc
\ifnum \tempcounta<\tempcountb \tempcounta=\tempcountb\fi
\divide \tempcountc by2
\horsize{\tempcountb}{#3}%
\advance \tempcountb by-\tempcountc
\ifnum \tempcountb>0 \advance \tempcountb by80\fi
\ifnum \tempcounta<\tempcountb \tempcounta=\tempcountb\fi
\advance\tempcounta by20
}

\def\putvector(#1,#2)(#3,#4)#5#6{{%
\xpos=#1
\ypos=#2
\run=#3
\rise=#4
\arrowlength=#5
\arrowtype=#6
\ifnum \arrowtype<0
    \ifnum \run=0
        \advance \ypos by-\arrowlength
    \else
        \tempcounta \arrowlength
        \multiply \tempcounta by\rise
        \divide \tempcounta by\run
        \ifnum\run>0
            \advance \xpos by\arrowlength
            \advance \ypos by\tempcounta
        \else
            \advance \xpos by-\arrowlength
            \advance \ypos by-\tempcounta
        \fi
    \fi
    \multiply \arrowtype by-1
    \multiply \rise by-1
    \multiply \run by-1
\fi
\ifnum \arrowtype=1
    \put(\xpos,\ypos){\vector(\run,\rise){\arrowlength}}%
\else\ifnum \arrowtype=2
    \put(\xpos,\ypos){\mvector(\run,\rise)\arrowlength}%
\else\ifnum\arrowtype=3
    \put(\xpos,\ypos){\evector(\run,\rise){\arrowlength}}%
\fi\fi\fi
}}

\def\putsplitvector(#1,#2)#3#4{%%
\xpos #1
\ypos #2
\arrowtype #4
\halflength #3
\arrowlength #3
\gap 140
\advance \halflength by-\gap
\divide \halflength by2
\ifnum \arrowtype=1
    \put(\xpos,\ypos){\line(0,-1){\halflength}}%
    \advance\ypos by-\halflength
    \advance\ypos by-\gap
    \put(\xpos,\ypos){\vector(0,-1){\halflength}}%
\else\ifnum \arrowtype=2
    \put(\xpos,\ypos){\line(0,-1)\halflength}%
    \put(\xpos,\ypos){\vector(0,-1)3}%
    \advance\ypos by-\halflength
    \advance\ypos by-\gap
    \put(\xpos,\ypos){\vector(0,-1){\halflength}}%
\else\ifnum\arrowtype=3
    \put(\xpos,\ypos){\line(0,-1)\halflength}%
    \advance\ypos by-\halflength
    \advance\ypos by-\gap
    \put(\xpos,\ypos){\evector(0,-1){\halflength}}%
\else\ifnum \arrowtype=-1
    \advance \ypos by-\arrowlength
    \put(\xpos,\ypos){\line(0,1){\halflength}}%
    \advance\ypos by\halflength
    \advance\ypos by\gap
    \put(\xpos,\ypos){\vector(0,1){\halflength}}%
\else\ifnum \arrowtype=-2
    \advance \ypos by-\arrowlength
    \put(\xpos,\ypos){\line(0,1)\halflength}%
    \put(\xpos,\ypos){\vector(0,1)3}%
    \advance\ypos by\halflength
    \advance\ypos by\gap
    \put(\xpos,\ypos){\vector(0,1){\halflength}}%
\else\ifnum\arrowtype=-3
    \advance \ypos by-\arrowlength
    \put(\xpos,\ypos){\line(0,1)\halflength}%
    \advance\ypos by\halflength
    \advance\ypos by\gap
    \put(\xpos,\ypos){\evector(0,1){\halflength}}%
\fi\fi\fi\fi\fi\fi
}

\def\putmorphism(#1)(#2,#3)[#4`#5`#6]#7#8#9{{%
\run #2
\rise #3
\ifnum\rise=0
  \puthmorphism(#1)[#4`#5`#6]{#7}{#8}{#9}%
\else\ifnum\run=0
  \putvmorphism(#1)[#4`#5`#6]{#7}{#8}{#9}%
\else
\setpos(#1)%
\arrowlength #7
\arrowtype #8
\ifnum\run=0
\else\ifnum\rise=0
\else
\ifnum\run>0
    \coefa=1
\else
   \coefa=-1
\fi
\ifnum\arrowtype>0
   \coefb=0
   \coefc=-1
\else
   \coefb=\coefa
   \coefc=1
   \arrowtype=-\arrowtype
\fi
\width=2
\multiply \width by\run
\divide \width by\rise
\ifnum \width<0  \width=-\width\fi
\advance\width by60
\if l#9 \width=-\width\fi
\putbox(\xpos,\ypos){#4}%            %node 1
{\multiply \coefa by\arrowlength%      %node 2
\advance\xpos by\coefa
\multiply \coefa by\rise
\divide \coefa by\run
\advance \ypos by\coefa
\putbox(\xpos,\ypos){#5} }%
{\multiply \coefa by\arrowlength%      %label
\divide \coefa by2
\advance \xpos by\coefa
\advance \xpos by\width
\multiply \coefa by\rise
\divide \coefa by\run
\advance \ypos by\coefa
\if l#9%
   \put(\xpos,\ypos){\makebox(0,0)[r]{$#6$}}%
\else\if r#9%
   \put(\xpos,\ypos){\makebox(0,0)[l]{$#6$}}%
\fi\fi }%
{\multiply \rise by-\coefc%             %arrow
\multiply \run by-\coefc
\multiply \coefb by\arrowlength
\advance \xpos by\coefb
\multiply \coefb by\rise
\divide \coefb by\run
\advance \ypos by\coefb
\multiply \coefc by70
\advance \ypos by\coefc
\multiply \coefc by\run
\divide \coefc by\rise
\advance \xpos by\coefc
\multiply \coefa by140
\multiply \coefa by\run
\divide \coefa by\rise
\advance \arrowlength by\coefa
\ifnum \arrowtype=1
   \put(\xpos,\ypos){\vector(\run,\rise){\arrowlength}}%
\else\ifnum\arrowtype=2
   \put(\xpos,\ypos){\mvector(\run,\rise){\arrowlength}}%
\else\ifnum\arrowtype=3
   \put(\xpos,\ypos){\evector(\run,\rise){\arrowlength}}%
\fi\fi\fi}\fi\fi\fi\fi}}

\def\puthmorphism(#1,#2)[#3`#4`#5]#6#7#8{{%
\xpos #1
\ypos #2
\width #6
\arrowlength #6
\putbox(\xpos,\ypos){#3\vphantom{#4}}%
{\advance \xpos by\arrowlength
\putbox(\xpos,\ypos){\vphantom{#3}#4}}%
\horsize{\tempcounta}{#3}%
\horsize{\tempcountb}{#4}%
\divide \tempcounta by2
\divide \tempcountb by2
\advance \tempcounta by30
\advance \tempcountb by30
\advance \xpos by\tempcounta
\advance \arrowlength by-\tempcounta
\advance \arrowlength by-\tempcountb
\putvector(\xpos,\ypos)(1,0){\arrowlength}{#7}%
\divide \arrowlength by2
\advance \xpos by\arrowlength
\vertsize{\tempcounta}{#5}%
\divide\tempcounta by2
\advance \tempcounta by20
\if a#8 %
   \advance \ypos by\tempcounta
   \putbox(\xpos,\ypos){#5}%
\else
   \advance \ypos by-\tempcounta
   \putbox(\xpos,\ypos){#5}%
\fi}}

\def\putvmorphism(#1,#2)[#3`#4`#5]#6#7#8{{%
\xpos #1
\ypos #2
\arrowlength #6
\arrowtype #7
\settowidth{\xlen}{$#5$}%
\putbox(\xpos,\ypos){#3}%
{\advance \ypos by-\arrowlength
\putbox(\xpos,\ypos){#4}}%
{\advance\arrowlength by-140
\advance \ypos by-70
\ifdim\xlen>0pt
   \if m#8%
      \putsplitvector(\xpos,\ypos){\arrowlength}{\arrowtype}%
   \else
      \putvector(\xpos,\ypos)(0,-1){\arrowlength}{\arrowtype}%
   \fi
\else
   \putvector(\xpos,\ypos)(0,-1){\arrowlength}{\arrowtype}%
\fi}%
\ifdim\xlen>0pt
   \divide \arrowlength by2
   \advance\ypos by-\arrowlength
   \if l#8%
      \advance \xpos by-40
      \put(\xpos,\ypos){\makebox(0,0)[r]{$#5$}}%
   \else\if r#8%
      \advance \xpos by40
      \put(\xpos,\ypos){\makebox(0,0)[l]{$#5$}}%
   \else
      \putbox(\xpos,\ypos){#5}%
   \fi\fi
\fi
}}

\def\topadjust[#1`#2`#3]{%
\yoff=10
\vertadjust[#1`#2`{#3}]%
\advance \yext by\tempcounta
\advance \yext by 10
}
\def\botadjust[#1`#2`#3]{%
\vertadjust[#1`#2`{#3}]%
\advance \yext by\tempcounta
\advance \yoff by-\tempcounta
}
\def\leftadjust[#1`#2`#3]{%
\xoff=0
\horadjust[#1`#2`{#3}]%
\advance \xext by\tempcounta
\advance \xoff by-\tempcounta
}
\def\rightadjust[#1`#2`#3]{%
\horadjust[#1`#2`{#3}]%
\advance \xext by\tempcounta
}
\def\rightsladjust[#1`#2`#3]{%
\sladjust[#1`#2`{#3}]{\width}%
\advance \xext by\tempcounta
}
\def\leftsladjust[#1`#2`#3]{%
\xoff=0
\sladjust[#1`#2`{#3}]{\width}%
\advance \xext by\tempcounta
\advance \xoff by-\tempcounta
}
\def\adjust[#1`#2;#3`#4;#5`#6;#7`#8]{%
\topadjust[#1``{#2}]
\leftadjust[#3``{#4}]
\rightadjust[#5``{#6}]
\botadjust[#7``{#8}]}

\def\putsquarep<#1>(#2)[#3;#4`#5`#6`#7]{{%
\setsqparms[#1]%
\setpos(#2)%
\settokens[#3]%
\puthmorphism(\xpos,\ypos)[\tokenc`\tokend`{#7}]{\width}{\arrowtyped}b%
\advance\ypos by \height
\puthmorphism(\xpos,\ypos)[\tokena`\tokenb`{#4}]{\width}{\arrowtypea}a%
\putvmorphism(\xpos,\ypos)[``{#5}]{\height}{\arrowtypeb}l%
\advance\xpos by \width
\putvmorphism(\xpos,\ypos)[``{#6}]{\height}{\arrowtypec}r%
}}

\def\putsquare{\@ifnextchar <{\putsquarep}{\putsquarep%
   <\arrowtypea`\arrowtypeb`\arrowtypec`\arrowtyped;\width`\height>}}
\def\square{\@ifnextchar< {\squarep}{\squarep
   <\arrowtypea`\arrowtypeb`\arrowtypec`\arrowtyped;\width`\height>}}
\def\squarep<#1>[#2`#3`#4`#5;#6`#7`#8`#9]{{%          %     #2------>#3
\setsqparms[#1]%                                      %      |       |
\xext=\width                                          %      |       |
\yext=\height                                         %    #7|       |#8
\topadjust[#2`#3`{#6}]%                               %      |       |
\botadjust[#4`#5`{#9}]%                               %      |       |
\leftadjust[#2`#4`{#7}]%                              %
\rightadjust[#3`#5`{#8}]%                             %     #4------>#5
\begin{picture}(\xext,\yext)(\xoff,\yoff)%                      #9
\putsquarep<\arrowtypea`\arrowtypeb`\arrowtypec`\arrowtyped;\width`\height>%
(0,0)[#2`#3`#4`#5;#6`#7`#8`{#9}]%
\end{picture}%
}}

\def\putptrianglep<#1>(#2,#3)[#4`#5`#6;#7`#8`#9]{{%
\settriparms[#1]%
\xpos=#2 \ypos=#3
\advance\ypos by \height
\puthmorphism(\xpos,\ypos)[#4`#5`{#7}]{\height}{\arrowtypea}a%
\putvmorphism(\xpos,\ypos)[`#6`{#8}]{\height}{\arrowtypeb}l%
\advance\xpos by\height
\putmorphism(\xpos,\ypos)(-1,-1)[``{#9}]{\height}{\arrowtypec}r%
}}

\def\putptriangle{\@ifnextchar <{\putptrianglep}{\putptrianglep
   <\arrowtypea`\arrowtypeb`\arrowtypec;\height>}}
\def\ptriangle{\@ifnextchar <{\ptrianglep}{\ptrianglep
   <\arrowtypea`\arrowtypeb`\arrowtypec;\height>}}

\def\ptrianglep<#1>[#2`#3`#4;#5`#6`#7]{{%%       #5
\settriparms[#1]%
\width=\height                         %      #2----->#3
\xext=\width                           %      |      /
\yext=\width                           %      |     /
\topadjust[#2`#3`{#5}]%                %    #6|    /#7
\botadjust[#3``]%                      %      |   /
\leftadjust[#2`#4`{#6}]%               %      |  /
\rightsladjust[#3`#4`{#7}]%            %
\begin{picture}(\xext,\yext)(\xoff,\yoff)%    #4
\putptrianglep<\arrowtypea`\arrowtypeb`\arrowtypec;\height>%
(0,0)[#2`#3`#4;#5`#6`{#7}]%
\end{picture}%
}}

\def\putqtrianglep<#1>(#2,#3)[#4`#5`#6;#7`#8`#9]{{%
\settriparms[#1]%
\xpos=#2 \ypos=#3
\advance\ypos by\height
\puthmorphism(\xpos,\ypos)[#4`#5`{#7}]{\height}{\arrowtypea}a%
\putmorphism(\xpos,\ypos)(1,-1)[``{#8}]{\height}{\arrowtypeb}l%
\advance\xpos by\height
\putvmorphism(\xpos,\ypos)[`#6`{#9}]{\height}{\arrowtypec}r%
}}

\def\putqtriangle{\@ifnextchar <{\putqtrianglep}{\putqtrianglep
   <\arrowtypea`\arrowtypeb`\arrowtypec;\height>}}
\def\qtriangle{\@ifnextchar <{\qtrianglep}{\qtrianglep
   <\arrowtypea`\arrowtypeb`\arrowtypec;\height>}}

\def\qtrianglep<#1>[#2`#3`#4;#5`#6`#7]{{%%
\settriparms[#1]%                                  #5
\width=\height                         %        #2----->#3
\xext=\width                           %         \      |
\yext=\height                          %          \     |
\topadjust[#2`#3`{#5}]%                %         #6\    |#7
\botadjust[#4``]%                      %            \   |
\leftsladjust[#2`#4`{#6}]%             %             \  |
\rightadjust[#3`#4`{#7}]%              %
\begin{picture}(\xext,\yext)(\xoff,\yoff)%             #4
\putqtrianglep<\arrowtypea`\arrowtypeb`\arrowtypec;\height>%
(0,0)[#2`#3`#4;#5`#6`{#7}]%
\end{picture}%
}}

\def\putdtrianglep<#1>(#2,#3)[#4`#5`#6;#7`#8`#9]{{%
\settriparms[#1]%
\xpos=#2 \ypos=#3
\puthmorphism(\xpos,\ypos)[#5`#6`{#9}]{\height}{\arrowtypec}b%
\advance\xpos by \height \advance\ypos by\height
\putmorphism(\xpos,\ypos)(-1,-1)[``{#7}]{\height}{\arrowtypea}l%
\putvmorphism(\xpos,\ypos)[#4``{#8}]{\height}{\arrowtypeb}r%
}}

\def\putdtriangle{\@ifnextchar <{\putdtrianglep}{\putdtrianglep
   <\arrowtypea`\arrowtypeb`\arrowtypec;\height>}}
\def\dtriangle{\@ifnextchar <{\dtrianglep}{\dtrianglep
   <\arrowtypea`\arrowtypeb`\arrowtypec;\height>}}

\def\dtrianglep<#1>[#2`#3`#4;#5`#6`#7]{{%%
\settriparms[#1]%                                          #2
\width=\height                         %                  / |
\xext=\width                           %                 /  |
\yext=\height                          %              #5/   |#6
\topadjust[#2``]%                      %               /    |
\botadjust[#3`#4`{#7}]%                %              /     |
\leftsladjust[#3`#2`{#5}]%             %
\rightadjust[#2`#4`{#6}]%              %            #3----->#4
\begin{picture}(\xext,\yext)(\xoff,\yoff)%              #7
\putdtrianglep<\arrowtypea`\arrowtypeb`\arrowtypec;\height>%
(0,0)[#2`#3`#4;#5`#6`{#7}]%
\end{picture}%
}}

\def\putbtrianglep<#1>(#2,#3)[#4`#5`#6;#7`#8`#9]{{%
\settriparms[#1]%
\xpos=#2 \ypos=#3
\puthmorphism(\xpos,\ypos)[#5`#6`{#9}]{\height}{\arrowtypec}b%
\advance\ypos by\height
\putmorphism(\xpos,\ypos)(1,-1)[``{#8}]{\height}{\arrowtypeb}r%
\putvmorphism(\xpos,\ypos)[#4``{#7}]{\height}{\arrowtypea}l%
}}

\def\putbtriangle{\@ifnextchar <{\putbtrianglep}{\putbtrianglep
   <\arrowtypea`\arrowtypeb`\arrowtypec;\height>}}
\def\btriangle{\@ifnextchar <{\btrianglep}{\btrianglep
   <\arrowtypea`\arrowtypeb`\arrowtypec;\height>}}

\def\btrianglep<#1>[#2`#3`#4;#5`#6`#7]{{%%
\settriparms[#1]%                                     #2
\width=\height                         %              | \
\xext=\width                           %              |  \
\yext=\height                          %            #5|   \#6
\topadjust[#2``]%                      %              |    \
\botadjust[#3`#4`{#7}]%                %              |     \
\leftadjust[#2`#3`{#5}]%               %
\rightsladjust[#4`#2`{#6}]%            %              #3----->#4
\begin{picture}(\xext,\yext)(\xoff,\yoff)%                #7
\putbtrianglep<\arrowtypea`\arrowtypeb`\arrowtypec;\height>%
(0,0)[#2`#3`#4;#5`#6`{#7}]%
\end{picture}%
}}

\def\putAtrianglep<#1>(#2,#3)[#4`#5`#6;#7`#8`#9]{{%
\settriparms[#1]%
\xpos=#2 \ypos=#3
{\multiply \height by2
\puthmorphism(\xpos,\ypos)[#5`#6`{#9}]{\height}{\arrowtypec}b}%
\advance\xpos by\height \advance\ypos by\height
\putmorphism(\xpos,\ypos)(-1,-1)[#4``{#7}]{\height}{\arrowtypea}l%
\putmorphism(\xpos,\ypos)(1,-1)[``{#8}]{\height}{\arrowtypeb}r%
}}

\def\putAtriangle{\@ifnextchar <{\putAtrianglep}{\putAtrianglep
   <\arrowtypea`\arrowtypeb`\arrowtypec;\height>}}
\def\Atriangle{\@ifnextchar <{\Atrianglep}{\Atrianglep
   <\arrowtypea`\arrowtypeb`\arrowtypec;\height>}}

\def\Atrianglep<#1>[#2`#3`#4;#5`#6`#7]{{%%
\settriparms[#1]%                                 #2
\width=\height                         %         /   \
\xext=\width                           %        /     \
\yext=\height                          %     #5/       \#6
\topadjust[#2``]%                      %      /         \
\botadjust[#3`#4`{#7}]%                %     /           \
\multiply \xext by2 %                  %
\leftsladjust[#3`#2`{#5}]%             %   #3------------>#4
\rightsladjust[#4`#2`{#6}]%            %          #7
\begin{picture}(\xext,\yext)(\xoff,\yoff)%
\putAtrianglep<\arrowtypea`\arrowtypeb`\arrowtypec;\height>%
(0,0)[#2`#3`#4;#5`#6`{#7}]%
\end{picture}%
}}

\def\putAtrianglepairp<#1>(#2)[#3;#4`#5`#6`#7`#8]{{
\settripairparms[#1]%
\setpos(#2)%
\settokens[#3]%
\puthmorphism(\xpos,\ypos)[\tokenb`\tokenc`{#7}]{\height}{\arrowtyped}b%
\advance\xpos by\height
\advance\ypos by\height
\putmorphism(\xpos,\ypos)(-1,-1)[\tokena``{#4}]{\height}{\arrowtypea}l%
\putvmorphism(\xpos,\ypos)[``{#5}]{\height}{\arrowtypeb}m%
\putmorphism(\xpos,\ypos)(1,-1)[``{#6}]{\height}{\arrowtypec}r%
}}

\def\putAtrianglepair{\@ifnextchar <{\putAtrianglepairp}{\putAtrianglepairp%
   <\arrowtypea`\arrowtypeb`\arrowtypec`\arrowtyped`\arrowtypee;\height>}}
\def\Atrianglepair{\@ifnextchar <{\Atrianglepairp}{\Atrianglepairp%
   <\arrowtypea`\arrowtypeb`\arrowtypec`\arrowtyped`\arrowtypee;\height>}}

\def\Atrianglepairp<#1>[#2;#3`#4`#5`#6`#7]{{%
\settripairparms[#1]%
\settokens[#2]%
\width=\height
\xext=\width
\yext=\height
\topadjust[\tokena``]%
\vertadjust[\tokenb`\tokenc`{#6}]%                      %  #2a
\tempcountd=\tempcounta                       %           / | \
\vertadjust[\tokenc`\tokend`{#7}]%            %          /  |  \
\ifnum\tempcounta<\tempcountd                 %       #3/  #4   \#5
\tempcounta=\tempcountd\fi                    %        /    |    \
\advance \yext by\tempcounta                  %       /     |     \
\advance \yoff by-\tempcounta                 %
\multiply \xext by2 %                         %     #2b---->#2c---->#2d
\leftsladjust[\tokenb`\tokena`{#3}]%          %         #6     #7
\rightsladjust[\tokend`\tokena`{#5}]%
\begin{picture}(\xext,\yext)(\xoff,\yoff)%
\putAtrianglepairp
<\arrowtypea`\arrowtypeb`\arrowtypec`\arrowtyped`\arrowtypee;\height>%
(0,0)[#2;#3`#4`#5`#6`{#7}]%
\end{picture}%
}}

\def\putVtrianglep<#1>(#2,#3)[#4`#5`#6;#7`#8`#9]{{%
\settriparms[#1]%
\xpos=#2 \ypos=#3
\advance\ypos by\height
{\multiply\height by2
\puthmorphism(\xpos,\ypos)[#4`#5`{#7}]{\height}{\arrowtypea}a}%
\putmorphism(\xpos,\ypos)(1,-1)[`#6`{#8}]{\height}{\arrowtypeb}l%
\advance\xpos by\height
\advance\xpos by\height
\putmorphism(\xpos,\ypos)(-1,-1)[``{#9}]{\height}{\arrowtypec}r%
}}

\def\putVtriangle{\@ifnextchar <{\putVtrianglep}{\putVtrianglep
   <\arrowtypea`\arrowtypeb`\arrowtypec;\height>}}
\def\Vtriangle{\@ifnextchar <{\Vtrianglep}{\Vtrianglep
   <\arrowtypea`\arrowtypeb`\arrowtypec;\height>}}

\def\Vtrianglep<#1>[#2`#3`#4;#5`#6`#7]{{%%
\settriparms[#1]%                                      #5
\width=\height                         %        #2------------->#3
\xext=\width                           %         \             /
\yext=\height                          %          \           /
\topadjust[#2`#3`{#5}]%                %         #6\         /#7
\botadjust[#4``]%                      %            \       /
\multiply \xext by2 %                  %             \     /
\leftsladjust[#2`#3`{#6}]%             %
\rightsladjust[#3`#4`{#7}]%            %               #4
\begin{picture}(\xext,\yext)(\xoff,\yoff)%
\putVtrianglep<\arrowtypea`\arrowtypeb`\arrowtypec;\height>%
(0,0)[#2`#3`#4;#5`#6`{#7}]%
\end{picture}%
}}

\def\putVtrianglepairp<#1>(#2)[#3;#4`#5`#6`#7`#8]{{
\settripairparms[#1]%
\setpos(#2)%
\settokens[#3]%
\advance\ypos by\height
\putmorphism(\xpos,\ypos)(1,-1)[`\tokend`{#6}]{\height}{\arrowtypec}l%
\puthmorphism(\xpos,\ypos)[\tokena`\tokenb`{#4}]{\height}{\arrowtypea}a%
\advance\xpos by\height
\putvmorphism(\xpos,\ypos)[``{#7}]{\height}{\arrowtyped}m%
\advance\xpos by\height
\putmorphism(\xpos,\ypos)(-1,-1)[``{#8}]{\height}{\arrowtypee}r%
}}

\def\putVtrianglepair{\@ifnextchar <{\putVtrianglepairp}{\putVtrianglepairp%
    <\arrowtypea`\arrowtypeb`\arrowtypec`\arrowtyped`\arrowtypee;\height>}}
\def\Vtrianglepair{\@ifnextchar <{\Vtrianglepairp}{\Vtrianglepairp%
    <\arrowtypea`\arrowtypeb`\arrowtypec`\arrowtyped`\arrowtypee;\height>}}

\def\Vtrianglepairp<#1>[#2;#3`#4`#5`#6`#7]{{%
\settripairparms[#1]%
\settokens[#2]%                            #3      #4
\xext=\height                  %        #2a---->#2b---->#2c
\width=\height                 %         \      |      /
\yext=\height                  %          \     |     /
\vertadjust[\tokena`\tokenb`{#4}]%       #5\   #6    /#7
\tempcountd=\tempcounta        %            \   |   /
\vertadjust[\tokenb`\tokenc`{#5}]%           \  |  /
\ifnum\tempcounta<\tempcountd%
\tempcounta=\tempcountd\fi%                    #2d
\advance \yext by\tempcounta
\botadjust[\tokend``]%
\multiply \xext by2
\leftsladjust[\tokena`\tokend`{#6}]%
\rightsladjust[\tokenc`\tokend`{#7}]%
\begin{picture}(\xext,\yext)(\xoff,\yoff)%
\putVtrianglepairp
<\arrowtypea`\arrowtypeb`\arrowtypec`\arrowtyped`\arrowtypee;\height>%
(0,0)[#2;#3`#4`#5`#6`{#7}]%
\end{picture}%
}}

\def\putCtrianglep<#1>(#2,#3)[#4`#5`#6;#7`#8`#9]{{%
\settriparms[#1]%
\xpos=#2 \ypos=#3
\advance\ypos by\height
\putmorphism(\xpos,\ypos)(1,-1)[``{#9}]{\height}{\arrowtypec}l%
\advance\xpos by\height
\advance\ypos by\height
\putmorphism(\xpos,\ypos)(-1,-1)[#4`#5`{#7}]{\height}{\arrowtypea}l%
{\multiply\height by 2
\putvmorphism(\xpos,\ypos)[`#6`{#8}]{\height}{\arrowtypeb}r}%
}}

\def\putCtriangle{\@ifnextchar <{\putCtrianglep}{\putCtrianglep
    <\arrowtypea`\arrowtypeb`\arrowtypec;\height>}}
\def\Ctriangle{\@ifnextchar <{\Ctrianglep}{\Ctrianglep
    <\arrowtypea`\arrowtypeb`\arrowtypec;\height>}}

\def\Ctrianglep<#1>[#2`#3`#4;#5`#6`#7]{{%%
\settriparms[#1]%                                         #2
\width=\height                          %                / |
\xext=\width                            %               /  |
\yext=\height                           %            #5/   |
\multiply \yext by2 %                   %             /    |
\topadjust[#2``]%                       %            /     |
\botadjust[#4``]%                       %           v      |
\sladjust[#3`#2`{#5}]{\width}%          %          #3      |#6
\tempcountd=\tempcounta                 %           \      |
\sladjust[#3`#4`{#7}]{\width}%          %            \     |
\ifnum \tempcounta<\tempcountd          %           #7\    |
\tempcounta=\tempcountd\fi              %              \   |
\advance \xext by\tempcounta            %               \  |
\advance \xoff by-\tempcounta           %
\rightadjust[#2`#4`{#6}]%               %                 #4
\begin{picture}(\xext,\yext)(\xoff,\yoff)%
\putCtrianglep<\arrowtypea`\arrowtypeb`\arrowtypec;\height>%
(0,0)[#2`#3`#4;#5`#6`{#7}]%
\end{picture}%
}}

\def\putDtrianglep<#1>(#2,#3)[#4`#5`#6;#7`#8`#9]{{%
\settriparms[#1]%
\xpos=#2 \ypos=#3
\advance\xpos by\height \advance\ypos by\height
\putmorphism(\xpos,\ypos)(-1,-1)[``{#9}]{\height}{\arrowtypec}r%
\advance\xpos by-\height \advance\ypos by\height
\putmorphism(\xpos,\ypos)(1,-1)[`#5`{#8}]{\height}{\arrowtypeb}r%
{\multiply\height by 2
\putvmorphism(\xpos,\ypos)[#4`#6`{#7}]{\height}{\arrowtypea}l}%
}}

\def\putDtriangle{\@ifnextchar <{\putDtrianglep}{\putDtrianglep
    <\arrowtypea`\arrowtypeb`\arrowtypec;\height>}}
\def\Dtriangle{\@ifnextchar <{\Dtrianglep}{\Dtrianglep
   <\arrowtypea`\arrowtypeb`\arrowtypec;\height>}}

\def\Dtrianglep<#1>[#2`#3`#4;#5`#6`#7]{{%%
\settriparms[#1]%                                 #2
\width=\height                         %          | \
\xext=\height                          %          |  \
\yext=\height                          %          |   \#6
\multiply \yext by2 %                  %          |    \
\topadjust[#2``]%                      %          |     \
\botadjust[#4``]%                      %          |
\leftadjust[#2`#4`{#5}]%               %        #5|      #3
\sladjust[#3`#2`{#5}]{\height}%        %          |      /
\tempcountd=\tempcountd                %          |     /
\sladjust[#3`#4`{#7}]{\height}%        %          |    /#7
\ifnum \tempcounta<\tempcountd         %          |   /
\tempcounta=\tempcountd\fi             %          |  /
\advance \xext by\tempcounta           %
\begin{picture}(\xext,\yext)(\xoff,\yoff)%        #4
\putDtrianglep<\arrowtypea`\arrowtypeb`\arrowtypec;\height>%
(0,0)[#2`#3`#4;#5`#6`{#7}]%
\end{picture}%
}}

\def\setrecparms[#1`#2]{\width=#1 \height=#2}%
\def\recursep<#1`#2>[#3;#4`#5`#6`#7`#8]{{%
\width=#1 \height=#2
\settokens[#3]
\settowidth{\tempdimen}{$\tokena$}
\ifdim\tempdimen=0pt
  \savebox{\tempboxa}{\hbox{$\tokenb$}}%
  \savebox{\tempboxb}{\hbox{$\tokend$}}%
  \savebox{\tempboxc}{\hbox{$#6$}}%
\else
  \savebox{\tempboxa}{\hbox{$\hbox{$\tokena$}\times\hbox{$\tokenb$}$}}%
  \savebox{\tempboxb}{\hbox{$\hbox{$\tokena$}\times\hbox{$\tokend$}$}}%
  \savebox{\tempboxc}{\hbox{$\hbox{$\tokena$}\times\hbox{$#6$}$}}%
\fi
\ypos=\height
\divide\ypos by 2
\xpos=\ypos
\advance\xpos by \width
\xext=\xpos \yext=\height
\topadjust[#3`\usebox{\tempboxa}`{#4}]%
\botadjust[#5`\usebox{\tempboxb}`{#8}]%
\sladjust[\tokenc`\tokenb`{#5}]{\ypos}%
\tempcountd=\tempcounta
\sladjust[\tokenc`\tokend`{#5}]{\ypos}%
\ifnum \tempcounta<\tempcountd
\tempcounta=\tempcountd\fi
\advance \xext by\tempcounta
\advance \xoff by-\tempcounta
\rightadjust[\usebox{\tempboxa}`\usebox{\tempboxb}`\usebox{\tempboxc}]%
\bfig
\putCtrianglep<-1`1`1;\ypos>(0,0)[`\tokenc`;#5`#6`{#7}]%
\puthmorphism(\ypos,0)[\tokend`\usebox{\tempboxb}`{#8}]{\width}{-1}b%
\puthmorphism(\ypos,\height)[\tokenb`\usebox{\tempboxa}`{#4}]{\width}{-1}a%
\advance\ypos by \width
\putvmorphism(\ypos,\height)[``\usebox{\tempboxc}]{\height}1r%
\efig
}}

\def\recurse{\@ifnextchar <{\recursep}{\recursep<\width`\height>}}

\def\puttwohmorphisms(#1,#2)[#3`#4;#5`#6]#7#8#9{{%
\puthmorphism(#1,#2)[#3`#4`]{#7}0a
\ypos=#2
\advance\ypos by 20
\puthmorphism(#1,\ypos)[\phantom{#3}`\phantom{#4}`#5]{#7}{#8}a
\advance\ypos by -40
\puthmorphism(#1,\ypos)[\phantom{#3}`\phantom{#4}`#6]{#7}{#9}b
}}

\def\puttwovmorphisms(#1,#2)[#3`#4;#5`#6]#7#8#9{{%
\putvmorphism(#1,#2)[#3`#4`]{#7}0a
\xpos=#1
\advance\xpos by -20
\putvmorphism(\xpos,#2)[\phantom{#3}`\phantom{#4}`#5]{#7}{#8}l
\advance\xpos by 40
\putvmorphism(\xpos,#2)[\phantom{#3}`\phantom{#4}`#6]{#7}{#9}r
}}

\def\puthcoequalizer(#1)[#2`#3`#4;#5`#6`#7]#8#9{{%
\setpos(#1)%
\puttwohmorphisms(\xpos,\ypos)[#2`#3;#5`#6]{#8}11%
\advance\xpos by #8
\puthmorphism(\xpos,\ypos)[\phantom{#3}`#4`#7]{#8}1{#9}
}}

\def\putvcoequalizer(#1)[#2`#3`#4;#5`#6`#7]#8#9{{%
\setpos(#1)%
\puttwovmorphisms(\xpos,\ypos)[#2`#3;#5`#6]{#8}11%
\advance\ypos by -#8
\putvmorphism(\xpos,\ypos)[\phantom{#3}`#4`#7]{#8}1{#9}
}}

\def\putthreehmorphisms(#1)[#2`#3;#4`#5`#6]#7(#8)#9{{%
\setpos(#1) \settypes(#8)
\if a#9 %
     \vertsize{\tempcounta}{#5}%
     \vertsize{\tempcountb}{#6}%
     \ifnum \tempcounta<\tempcountb \tempcounta=\tempcountb \fi
\else
     \vertsize{\tempcounta}{#4}%
     \vertsize{\tempcountb}{#5}%
     \ifnum \tempcounta<\tempcountb \tempcounta=\tempcountb \fi
\fi
\advance \tempcounta by 60
\puthmorphism(\xpos,\ypos)[#2`#3`#5]{#7}{\arrowtypeb}{#9}
\advance\ypos by \tempcounta
\puthmorphism(\xpos,\ypos)[\phantom{#2}`\phantom{#3}`#4]{#7}{\arrowtypea}{#9}
\advance\ypos by -\tempcounta \advance\ypos by -\tempcounta
\puthmorphism(\xpos,\ypos)[\phantom{#2}`\phantom{#3}`#6]{#7}{\arrowtypec}{#9}
}}

\def\putarc(#1,#2)[#3`#4`#5]#6#7#8{{%
\xpos #1
\ypos #2
\width #6
\arrowlength #6
\putbox(\xpos,\ypos){#3\vphantom{#4}}%
{\advance \xpos by\arrowlength
\putbox(\xpos,\ypos){\vphantom{#3}#4}}%
\horsize{\tempcounta}{#3}%
\horsize{\tempcountb}{#4}%
\divide \tempcounta by2
\divide \tempcountb by2
\advance \tempcounta by30
\advance \tempcountb by30
\advance \xpos by\tempcounta
\advance \arrowlength by-\tempcounta
\advance \arrowlength by-\tempcountb
\halflength=\arrowlength \divide\halflength by 2
\divide\arrowlength by 5
\put(\xpos,\ypos){\bezier{\arrowlength}(0,0)(50,50)(\halflength,50)}
\ifnum #7=-1 \put(\xpos,\ypos){\vector(-3,-2)0} \fi
\advance\xpos by \halflength
\put(\xpos,\ypos){\xpos=\halflength \advance\xpos by -50
   \bezier{\arrowlength}(0,50)(\xpos,50)(\halflength,0)}
\ifnum #7=1 {\advance \xpos by
   \halflength \put(\xpos,\ypos){\vector(3,-2)0}} \fi
\advance\ypos by 50
\vertsize{\tempcounta}{#5}%
\divide\tempcounta by2
\advance \tempcounta by20
\if a#8 %
   \advance \ypos by\tempcounta
   \putbox(\xpos,\ypos){#5}%
\else
   \advance \ypos by-\tempcounta
   \putbox(\xpos,\ypos){#5}%
\fi
}}

\makeatother

\usepackage{amsthm}
\usepackage{dsfont}
\usepackage{stmaryrd}
\newtheorem{theorem}{Theorem}[section]
\newtheorem{lemma}[theorem]{Lemma}
\newtheorem{corollary}[theorem]{Corollary}

\newtheorem{proposition}[theorem]{Proposition}

\makeindex \makeglossary
\begin{document}

\sloppy

%commands
\newcommand{\nl}{\hspace{2cm}\\ }

\def\nec{\Box}
\def\pos{\Diamond}
\def\diam{{\tiny\Diamond}}

\def\lc{\lceil}
\def\rc{\rceil}
\def\lf{\lfloor}
\def\rf{\rfloor}
\def\lk{\langle}
\def\rk{\rangle}
\def\blk{\dot{\langle\!\!\langle}}
\def\brk{\dot{\rangle\!\!\rangle}}

\newcommand{\pa}{\parallel}
\newcommand{\lra}{\longrightarrow}
\newcommand{\hra}{\hookrightarrow}
\newcommand{\hla}{\hookleftarrow}
\newcommand{\ra}{\rightarrow}
\newcommand{\la}{\leftarrow}
\newcommand{\lla}{\longleftarrow}
\newcommand{\da}{\downarrow}
\newcommand{\ua}{\uparrow}
\newcommand{\dA}{\downarrow\!\!\!^\bullet}
\newcommand{\uA}{\uparrow\!\!\!_\bullet}
\newcommand{\Da}{\Downarrow}
\newcommand{\DA}{\Downarrow\!\!\!^\bullet}
\newcommand{\UA}{\Uparrow\!\!\!_\bullet}
\newcommand{\Ua}{\Uparrow}
\newcommand{\Lra}{\Longrightarrow}
\newcommand{\Ra}{\Rightarrow}
\newcommand{\Lla}{\Longleftarrow}
\newcommand{\La}{\Leftarrow}
\newcommand{\nperp}{\perp\!\!\!\!\!\setminus\;\;}
\newcommand{\pq}{\preceq}

\newcommand{\lms}{\longmapsto}
\newcommand{\ms}{\mapsto}
\newcommand{\subseteqnot}{\subseteq\hskip-4 mm_\not\hskip3 mm}

\def\o{{\omega}}

\def\bA{{\bf A}}
\def\bEM{{EM}}
\def\bM{{\bf M}}
\def\bN{{\bf N}}
\def\bC{{\bf C}}
\def\bI{{\bf I}}
\def\bK{{K}}
\def\bL{{\bf L}}
\def\bT{{\bf T}}
\def\bS{{\bf S}}
\def\bD{{\bf D}}
\def\bB{{\bf B}}
\def\bW{{\bf W}}
\def\bP{{\bf P}}
\def\bX{{\bf X}}
\def\bY{{\bf Y}}
\def\ba{{\bf a}}
\def\bb{{\bf b}}
\def\bc{{\bf c}}
\def\bd{{\bf d}}
\def\bh{{\bf h}}
\def\bi{{\bf i}}
\def\bj{{\bf j}}
\def\bk{{\bf k}}
\def\bm{{\bf m}}
\def\bn{{\bf n}}
\def\bp{{\bf p}}
\def\bq{{\bf q}}
\def\be{{\bf e}}
\def\br{{\bf r}}
\def\bi{{\bf i}}
\def\bs{{\bf s}}
\def\bt{{\bf t}}
\def\jeden{{\bf 1}}
\def\dwa{{\bf 2}}
\def\trzy{{\bf 3}}

\def\cB{{\cal B}}
\def\cA{{\cal A}}
\def\cC{{\cal C}}
\def\cD{{\cal D}}
\def\cE{{\cal E}}
\def\cEM{{\cal EM}}
\def\cF{{\cal F}}
\def\cG{{\cal G}}
\def\cI{{\cal I}}
\def\cJ{{\cal J}}
\def\cK{{\cal K}}
\def\cL{{\cal L}}
\def\cN{{\cal N}}
\def\cM{{\cal M}}
\def\cO{{\cal O}}
\def\cP{{\cal P}}
\def\cQ{{\cal Q}}
\def\cR{{\cal R}}
\def\cS{{\cal S}}
\def\cT{{\cal T}}
\def\cU{{\cal U}}
\def\cV{{\cal V}}
\def\cW{{\cal W}}
\def\cX{{\cal X}}
\def\cY{{\cal Y}}

%categories

%of functors and monads
\def\Mnd{{\bf Mnd}}
\def\AMnd{{\bf AnMnd}}
\def\An{{\bf An}}
\def\San{{\bf San}}
\def\PMnd{{\bf PolyMnd}}
\def\SanMnd{{\bf SanMnd}}
\def\RiMnd{{\bf RiMnd}}
\def\End{{\bf End}}

%of theories
\def\ET{\bf ET}
\def\RegET{\bf RegET}
\def\RET{\bf RegET}
\def\LrET{\bf LrET}
\def\RiET{\bf RiET}
\def\SregET{\bf SregET}

%of Lawvere theories
\def\LT{\bf LT}
\def\RegLT{\bf RegLT}
\def\ALT{\bf AnLT}
\def\RiLT{\bf RiLT}

%of Operads
\def\FOp{\bf FOp}
\def\RegOp{\bf RegOp}
\def\SOp{\bf SOp}
\def\RiOp{\bf RiOp}

%various other categories and such
\def\bCat{{{\bf Cat}}}
\def\MonCat{{{\bf MonCat}}}
\def\Mon{{{\bf Mon}}}
\def\Cat{{{\bf Cat}}}

\def\F{\mathds{F}}
\def\S{\mathds{S}}
\def\I{\mathds{I}}
\def\B{\mathds{B}}

%functors
\def\V{\mathds{V}}
\def\W{\mathds{W}}
\def\M{\mathds{M}}
\def\N{\mathds{N}}

\def\Op{{\cal O}p}

\def\Vb{\bar{\mathds{V}}}
\def\Wb{\bar{\mathds{W}}}

\def\Sym{{\cal S}ym}

%leftovers
\def\P{{\cal P}}
\def\Q{{\cal Q}}

\pagenumbering{arabic} \setcounter{page}{1}

\title{\bf\Large Theories of analytic monads}
\author{ Stanis\l aw Szawiel,\\ Marek Zawadowski}

\maketitle
\begin{abstract}  We characterize the equational theories and Lawvere theories that correspond to the categories of analytic and polynomial monads on $Set$, and hence also the categories of the symmetric and rigid operads in $Set$. We show that the category of analytic monads is equivalent to the category of regular-linear theories.  The category of polynomial monads is equivalent to the category of rigid theories, i.e. regular-linear theories satisfying an additional global condition. This solves a problem stated in \cite{CJ2}. The Lawvere theories corresponding to these monads are identified via some factorization systems.

\end{abstract}


%\newpage
\tableofcontents
\section{Introduction}

The category of algebras of a (finitary) equational theory can be equivalently described as a category of models of a Lawvere theory or as a category of algebras of a finitary monad on the category $Set$. In some cases there are also two other descriptions available. Some categories of algebras can be described also as algebras for a symmetric operad and some  can be described as algebras for a rigid\footnote{We call a `rigid operad' what was earlier called an `operad with non-standard amalgamations' (cf. \cite{HMP}). We decided to change the name since it is a very important notion deserving a handy name. The choice of the name is motivated by the property of the equational theories that correspond to such operads. The category of rigid operads can be identified with the full subcategory of the symmetric operad such that the actions of symmetric groups on their operations are free.} operad (c.f. [HMP], [LMF]). It is well known that the categories of equational theories $\ET$, Lawvere theories $\LT$ and monads (on $Set$) $\Mnd$ are equivalent\footnote{These correspondences preserves the notion of a `model', i.e. the corresponding categories of algebras in the suitable sense are canonically equivalent.}. Its is also known that the categories of symmetric and rigid operads are equivalent to the categories of analytic and polynomial monads, respectively; see \cite{Z1}. In this paper we give a description of the subcategories of  $\ET$ and of $\LT$ that correspond to the categories of symmetric and rigid operads.

The equational theories corresponding to analytic monads are linear-regular theories. A linear-regular theory is an equational theory that can be axiomatized by equations having the same variables on both sides, each variable occurring exactly once. A linear-regular theory $T$ is rigid iff whenever a linear-regular equation
\[ t(x_1,\ldots x_n)=t(x_{\sigma(1)},\ldots x_{\sigma(n)})\]
 is provable in $T$ then the permutation $\sigma$ is the identity permutation. In the above equation $t(x_1,\ldots x_n)$ denotes any term with $n$ different variables $x_1,\ldots x_n$, each one occurring exactly once and $t(x_{\sigma(1)},\ldots x_{\sigma(n)})$ denotes the same term $t$ but with variables permuted according to $\sigma$. For example, the theory of commutative monoids is not rigid as it contains the equation
\[ m(x_1,x_2)=m(x_2,x_1) \]
The category of polynomial monads $\PMnd$ corresponds to the category of rigid theories $\RiET$. The notion of a linear-regular theory was considered in universal algebra but the notion of a rigid theory as well as that of a linear-regular interpretation seems to be new.  If all the axioms of an equational theory are linear-regular, then the theory is linear-regular. However, the problem whether a finite set of linear-regular equations defines a rigid theory is undecidable, (cf. \cite{BSZ}).

We also give a characterization of the categories of Lawvere theories that correspond to the categories of analytic and polynomial monads. The category $\F^{op}$,  opposite of the skeleton of the category of finite sets is the initial Lawvere theory. Thus it has a unique morphism into any other Lawvere theory $\pi :\F^{op}\ra \bT$. The class of morphisms in the image of $\pi$ closed under isomorphisms is called the class of structural morphisms in $\bT$. The class of right orthogonal morphisms to the structural morphisms is the class of analytic morphisms in $\bT$. A Lawvere theory $\bT$ is analytic  iff the classes of structural and analytic  morphisms form a factorization system and the automorphisms of any object $n$ in $\bT$ are determined by the automorphisms of $1$. A Lawvere theory is rigid if it is analytic and the symmetric group actions act freely on analytic operations. We show that the categories $\ALT$ of analytic and $\RiLT$ of rigid Lawvere theories correspond to the categories of analytic and of polynomial monads.

The following diagram illustrates the relations between the categories mentioned above. The vertical lines denote adjoint equivalences. Thus up to equivalence of categories there are only three categories in it, one on each level. One equivalent to the category of all finitary monads on $Set$, second  equivalent to the category of analytic monads on $Set$, and third equivalent to the category of polynomial monads on $Set$.  These levels are denoted by letters $f$, $a$, and $p$, respectively. Thus all four columns of equational theories, Lawvere theories, monads and operads\footnote{The column for operads is a bit shorter in this paper but it can be extended as we will show in \cite{SZ2}.} are `level-wise' equivalent. These columns are denoted by letters $e$, $l$, $m$ and $o$, respectively. The vertical functors heading up are inclusions of subcategories. The lower functors are full inclusions and the upper are inclusions that are full on isomorphisms.  The vertical functors heading down, the right adjoints to those heading up, are monadic. All the squares in the diagram commute up to canonical isomorphisms.
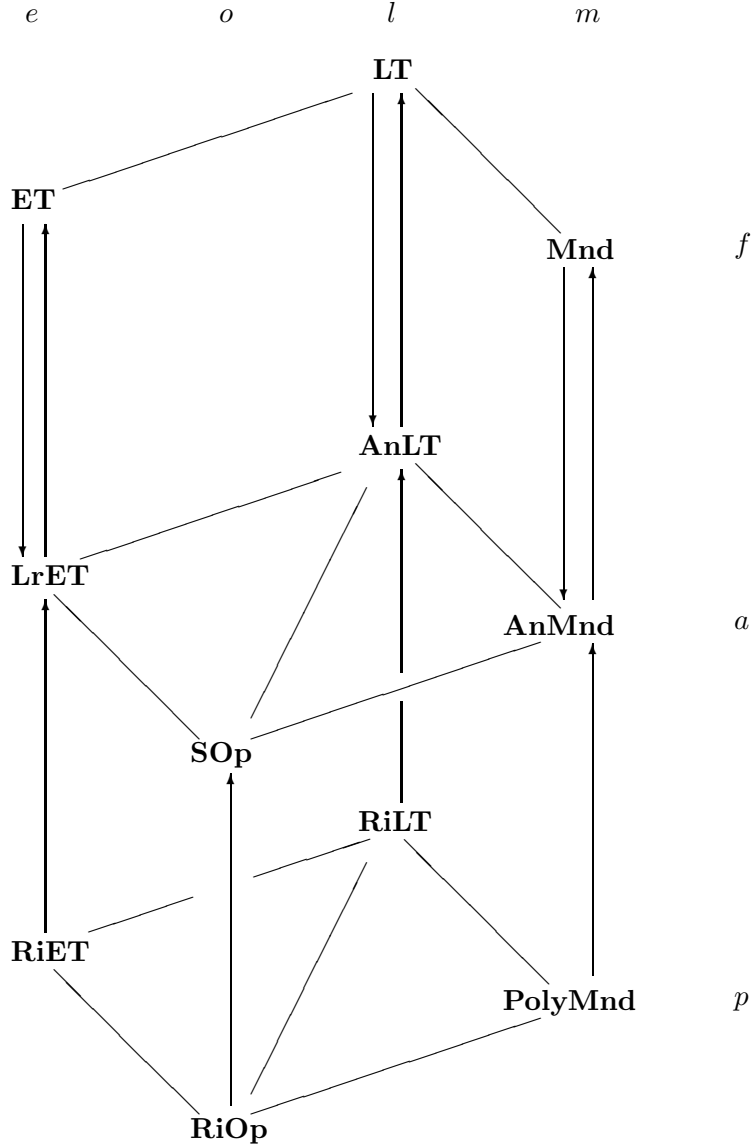
\begin{center} \xext=2800 \yext=4020
\begin{picture}(\xext,\yext)(\xoff,\yoff)
%level p
  %cats
  \put(620,40){$\RiOp$}
  \put(1700,480){$\PMnd$}
  \put(0,650){$\RiET$}
  \put(1200,1100){$\RiLT$}
   \put(2500,500){$p$}

 %horizontal
 \put(830,120){\line(3,1){1000}}
 \put(840,940){\line(3,1){400}}  \put(630,870){\line(-3,-1){360}}
 \put(650,120){\line(-1,1){500}}
 \put(1860,570){\line(-1,1){500}}
  %\put(790,150){\line(2,3){550}}
  \put(830,190){\line(1,2){400}}

 %vertical
   \put(120,750){\vector(0,1){1150}}
   \put(760,150){\vector(0,1){1150}}
   \put(2010,600){\vector(0,1){1150}}
    \put(1350,1200){\line(0,1){350}} \put(1350,1650){\vector(0,1){700}}

%level a
 %cats
  \put(620,1340){$\SOp$}
  \put(1700,1780){$\AMnd$}
  \put(0,1950){$\LrET$}
  \put(1200,2400){$\ALT$}
   \put(2500,1800){$a$}

 %horizontal
 \put(830,1420){\line(3,1){1000}}
 \put(240,2040){\line(3,1){900}}
 \put(650,1420){\line(-1,1){500}}
 \put(1900,1870){\line(-1,1){500}}
  \put(830,1490){\line(1,2){400}}

   %vertical
   \put(40,3200){\vector(0,-1){1150}} \put(120,2050){\vector(0,1){1150}}
   \put(1910,3050){\vector(0,-1){1150}} \put(2010,1900){\vector(0,1){1150}}

    \put(1250,3650){\vector(0,-1){1150}}
    \put(1350,2500){\vector(0,1){1150}}

%level f
 %cats
  \put(1850,3080){$\Mnd$}
  \put(0,3250){$\ET$}
  \put(1250,3700){$\LT$}
  \put(2500,3100){$f$}

 %horizontal
\put(180,3320){\line(3,1){1000}}
\put(1900,3170){\line(-1,1){500}}

\put(50,3900){$e$}
\put(720,3900){$o$}
\put(1300,3900){$l$}
\put(1950,3900){$m$}
\end{picture}
\end{center}
The notation concerning categories involved is displayed in the above diagram. The notation concerning functors is not on the diagram but it is meant to be systematically referring to the levels and columns they `connect'. The horizontal functors are denoted using letters from both columns they connect; the codomain by the script letter, the domain by its subscript, and the level is denoted by superscript. Thus the functor $\AMnd\ra \ALT$ will be denoted by $\cL_m^a$. We usually drop superscripts and often subscripts when it does not lead to confusion. Thus we can write, for example, $\cE=\cE_o=\cE_o^p : \RiOp \ra \RiET$. The vertical functors heading up are denoted by script letter $\cP$ with subscript indicating the column and superscript indicating the level of the codomain. The  vertical functors heading down are denoted by script letter $\cQ$ with subscript and superscript as those heading up. Thus we have, for example, functors $\cP=\cP^o=\cP^o_a: \RiOp \ra \SOp$ and $\cQ=\cQ_f=\cQ_f^m: \Mnd \ra \AMnd$. We will also refer to various diagonal morphisms and then we need to extend the notation concerning vertical functors by specifying both the columns of the domain and the codomain. For example, we write $\cP^{ol}_a: \SOp\ra \LT$ to denote one such functor and its right adjoint will be denoted by $\cQ^{lo}_f: \LT\ra \SOp$. In principle this notation will leave the codomain not always uniquely specified but in practice it it sufficient, and in fact usually much less is needed and each time it is used it will be recalled on the spot.

The paper is organized as follows. In Section 2 we recall categories of equational theories, Lawvere theories, monads on $Set$, and operads. We also discuss some of their subcategories. In Sections 3 we recall the correspondence between equational theories, Lawvere theories, and monads. In section 4 we study relations between Lawvere theories and operads.  We define a functor $\cL_o:\SOp \ra\LT$ from the category of symmetric operads the category of Lawvere theories. We identify its image and we show that its right adjoint is monadic. We also identify the image of the category of the rigid operads $\RiOp$ in $\LT$. In section 5 we relate the result from section 4 to monads. We note that finitary monads are monadic over analytic ones. But we also explain that this is a consequence of an even more fundamental fact that there is a lax monoidal monad on the category of analytic functors. This monoidal monad induces a distributive law that resembles what could be the formalization of the `combing trees' described in \cite{BD}. From this we obtain that finitary monads are monadic over analytic functors. This extends a result from \cite{Barr}.  In section 6 we define the embedding $\SOp$ in $\ET$ and characterize the images of both $\SOp$ and $\RiOp$. This gives the characterizations described at the beginning of the introduction that solves a problem stated in \cite{CJ2}. The paper ends with a section where we give some examples.

\subsection*{Notation}
We introduce notation that will be used in the whole paper.  $\o$  denotes the set of natural numbers. For $n\in\o$, we have $n=\{0,\ldots, n-1\}$, $[n]=\{0,\ldots, n\}$, $(n]=\{1,\ldots, n\}$. The set $X^n$ is interpreted as $X^{(n]}$ and it has a (natural) right action of the permutations group $S_n$ by composition. The skeletal category equivalent to the category of finite sets will be denoted by $\mathds{F}$. The objects of $\mathds{F}$ are sets $(n]$, for $n\in \o$.  The subcategories of $\mathds{F}$ with the same objects as $\mathds{F}$ but having as morphisms bijections, surjections, and injections will be denoted by $\mathds{B}$, $\mathds{S}$, $\mathds{I}$, respectively. When $S_n$ acts on the set $A$ on the right and on the set $B$ on the left, the set $A\otimes_nB$ is the usual tensor product of $S_n$-sets.

\section{Presentations of categories of algebras}

In this section we collect several categories whose objects describe (some) categories of algebras of finitary equational theories and whose morphisms
induce functors between such categories of algebras.

\subsection*{Equational theories}

By an {\em equational theory} we mean a pair of sets $T=(L,A)$, $L=\bigcup_{n\in \o} L_n$ and $L_n$ is the set of $n$-ary operations of $T$. The sets of operations of different arities are disjoint. The set $\cT r(L,\vec{x}^n)$ of terms of $L$ in context $\vec{x}^n=\lk x_1, \ldots, x_n \rk $ is the usual set of terms over $L$ build with the help of variables from $\vec{x}^n$. We write $t:\vec{x}^n$ for the term $t$ in context $\vec{x}^n$. Thus all the variables occurring in $t$ are among those in $\vec{x}^n$. The set $A$ is a set of equations in context $t=s : \vec{x}^n$, i.e. both $t:\vec{x}^n$ and $s:\vec{x}^n$ are terms in context.

A morphism of equational theories, an {\em interpretation},  $I : (L,A) \ra (L',A')$ is given by a set of functions $I_n : L_n \ra \cT r(L',\vec{x}^n)$, for $n\in \o$.  $I_n$'s extend to functions  $\bar{I}_n : \cT r(L,\vec{x}^n) \ra \cT r(L',\vec{x}^n)$, for $n\in \o$ as follows. We drop index $n$ in $\bar{I}_n$ when it does not lead to confusion.

\[ \bar{I}(x_i:\vec{x}^n) = x_i:\vec{x}^n \] for $i=1, \ldots, n$ and

\[ \bar{I}(f(t_1,\ldots, t_k) :\vec{x}^n) = I(f)(x_1\setminus \bar{I}(t_1),\ldots, x_k\setminus \bar{I}(t_k)):\vec{x}^n \]
for $f\in L_k$ and $t_i\in \cT r(L,\vec{x}^n)$ for $i=1,\ldots, k$. On the right-hand side we have a simultaneous substitution of terms $t_i$'s for variables $x_i$'s.  Moreover, for $I$ to be an interpretation we require that the equations are preserved, i.e. for any $t=s : \vec{x}^n$ in $A$ we have
\[ A' \vdash \bar{I}(t)=\bar{I}(s) : \vec{x}^n\]
where $A'\vdash$ is the provability in the equational logic from axioms in the set $A'$. We identify two such interpretations $I$ and  $I' : (L,A) \ra (L',A')$
iff they interpret all function symbols as provably equivalent terms, i.e.
\[ A' \vdash I(f)=I'(f) : \vec{x}^n\]
for any $n\in\o$ and $f\in L_n$. In this way we have defined a category of equational theories $\ET$.

A term in context $t:\vec{x}^n$ is {\em regular} if every variable in $\vec{x}^n$ occurs in $t$ at least once. A term in context $t:\vec{x}^n$ is {\em linear} if every variable in $\vec{x}^n$ occurs in $t$ at most once. A term in context $t:\vec{x}^n$ is {\em linear-regular} if it is both linear and regular.
An equation $s=t:\vec{x}^n$ is {\em regular} ({\em linear-regular}) iff both $s:\vec{x}^n$ and $t:\vec{x}^n$ are regular (linear-regular) terms in contexts.

A {\em simple $\phi$-substitution of a term in context}  $t:\vec{x}^n$ along a function $\phi:(n]\ra (k]$ is a term in context denoted $\phi\cdot t : \vec{x}^k$ such that every occurrence of the variable $x_i$ is replaced by the occurrence of $x_{\phi(i)}$. An {\em $\alpha$-conversion of a term in context}  $t:\vec{x}^n$ is a simple $\phi$-substitution of a term in context along a monomorphism $\phi:(n]\ra (k]$.

An equational theory $T=(L,A)$ is a {\em regular} ({\em linear-regular}) {\em theory} iff every equation $s=t:\vec{x}^n$ that is a consequence of the theory $T$ is a consequence of the set of regular (linear-regular) consequences of $T$. An interpretation is a {\em regular (linear-regular) interpretation} iff it interprets function symbols as regular (linear-regular) terms.

A theory $T=(L,A)$ is a {\em rigid theory} iff it is linear-regular and for any linear-regular term in context $t:\vec{x}^n$ whenever  $A\vdash t=\tau \cdot t:\vec{x}^n$ then $\tau$ is the identity permutation. $\tau \cdot t$ is the simple $\tau$-substitution of a term in context  $t:\vec{x}^n$ along a permutation $\tau:(n]\ra (n]\in S_n$.

Note that it is not assumed that the axioms of linear-regular theories are linear-regular. This is to keep the notion invariant under isomorphism of theories.
In particular, if $T=(L,A)$ is a linear-regular theory and $A'$ is the set of all equational consequences of the axioms from $A$ in the language $L$ then the $T'=(L,A')$ is also linear-regular. Of course $T'$ is isomorphic to $T$.  On the other hand, if the theory has linear-regular axioms then it is linear-regular.
Thus if we find a linear-regular set of axioms of an equational theory $T$ we can be sure that $T$ is linear-regular. However, it is not so easy to decide whether a theory is rigid. In fact, even if we have a finite linear-regular presentation of a theory it is still undecidable whether this theory is rigid (cf. \cite{BSZ}).

\vskip 2mm

{\em Remark.} We could also consider here strongly regular theories (c.f. \cite{CJ1}) that correspond to monotone monads (cf. \cite{Z1}). They are more specific than rigid theories. However this part of correspondence is of a bit different kind. The monotone monads are not just monads with certain additional properties but also with certain additional structure. The forgetful functor from monotone monads to monads (on $Set$) is not full on isomorphisms. These theories and some other theories of this kind will be treated elsewhere. We give the definition just to show the difference between the notions in the examples below.  A term in context $t:\vec{x}^n$ is a {\em strongly regular} iff it is linear-regular and the variables in the term $t$ occur in the same order as in the sequence  $\vec{x}^n$. An equation is $s=t:\vec{x}^n$ is a {\em strongly regular equation} iff both terms $s:\vec{x}^n$ and $t:\vec{x}^n$ are strongly regular. An equational theory $T=(L,A)$ is a {\em strongly regular theory} iff every equation $s=t:\vec{x}^n$ that is a consequence of the theory $T$ is a consequence of the set of strongly regular consequences of $T$. An interpretation is a {\em strongly regular interpretation} iff it interprets $n$-ary function symbols as strongly regular terms. One can easily see that  any strongly regular theory is rigid but the `embedding' functor $\SregET\lra \RiET$ is not even full on isomorphisms, where $\SregET$ denotes the category of strongly regular theories and strongly regular interpretations. The examples below show that there are rigid theories that are not strongly regular.

\vskip 2mm

We denote by $\LrET$ the subcategory of $\ET$ consisting of linear-regular theories and linear-regular interpretations. $\RiET$ denotes the full subcategory of $\LrET$ whose objects are rigid theories. $\RegET$ is a category of regular theories and regular interpretations. We have three inclusion functors
\[  \RiET \lra \LrET \lra \RegET \lra \ET \]
with the first inclusion being full and the other two being full on isomorphisms (cf. \cite{Z1}).

\vskip 2mm
{\em Examples.}
\begin{enumerate}
  \item The theory of monoids has two operations $m$ and $e$, of arity $2$ and $0$, respectively and equations
  \[ m(x_1,m(x_2,x_3))= m(m(x_1,x_2),x_3), \;\;\;\; m(x_1,e)=x_1=m(e,x_1) \]
  By the form of these equations, this theory is strongly regular and hence it is rigid as well.
   \item The theory of monoids with anti-involution has an additional unary operation $s$ and additional two axioms
  \[ m(s(x_1), s(x_2)) = s(m(x_2,x_1)), \;\;\;\; s(s(x_1))=x_1 \]
  This theory is not strongly regular but it can be shown that it is rigid.
  \item The theory of commutative monoids is the theory of monoids with an additional axiom
  \[ m(x_1,x_2) = m(x_2,x_1) \]
  Thus is it linear-regular by the form of the axioms but it is obviously not rigid.
  \item The theory of sup-lattices has two operations $\vee$ and $\bot$, of arity $2$ and $0$, respectively and equations
  \[ x_1\vee(x_2\vee x_3)= (x_1\vee x_2)\vee x_3,\;\;\;\; x_1\vee e=x_1= e \vee x_1 \]
  \[ x_1\vee x_1=x_1,  \;\;\;\; x_1\vee x_2= x_2 \vee x_1 \]
  This theory is regular but not linear.
  \item The theory of groups is not regular.
\end{enumerate}

\subsection*{Lawvere theories}
By a {\em Lawvere theory}, (cf. \cite{Lw}, \cite{KR}), we mean a category whose objects are natural numbers $\o$, so that $n$ is a product $1^n$  with chosen projections $\pi^n_i : n\ra 1$, for $n\in \o$ and $i\in (n]$. An interpretation (or a morphism) of Lawvere theories is a functor constant on objects, preserving the chosen projections.  Lawvere theories and their morphisms form a category that is denoted by $\LT$.

The initial object in the category $\LT$ is the category $\mathds{F}^{op}$ with the obvious inclusions as projections, see introduction. The unique morphism from $\mathds{F}^{op}$ into any Lawvere theory $\bT$ will be denoted by  $\pi: \mathds{F}^{op}\lra \bT$. Thus for $\phi: (n] \ra (m]$ in $\mathds{F}$  we have \mbox{$\pi_\phi= \lk \pi_{\phi(i)} \rk_{i\in (n]} : m \ra n$} in $\bT$.

Every equational theory has a model, and hence there is no inconsistent equational theory. But there are two equational theories that are nearly so. The terminal Lawvere theory $\mathds{1}$ has exactly one morphism between any two objects.  It has unique (up to isomorphism) one-element model. $\mathds{1}$ is not a regular theory. It also has some equivalent internal characterizations as the Lawvere theory (unique up to an isomorphism) which is a groupoid or where $0\cong 1$. The functor  $\pi: \mathds{F}^{op}\lra \mathds{1}$ is not faithful. There is yet another Lawvere theory with this property. It is a subtheory of $\mathds{1}$ in which there is no morphism $0\ra 1$. These categories are the only two Lawvere theories in which $2$ is the initial object.

The class of {\em structural morphisms} in $\bT$ is the closure under isomorphism of the image under $\pi$ of all morphisms in $\mathds{F}$. A morphism in $\bT$ is {\em analytic} iff it is right orthogonal to all structural morphisms.

By a factorization system in a category $\cC$ we mean the factorization system in the sense of \cite{FK}, see \cite{CJKP} sec 2.8, i.e. it consists of two classes of morphisms in $\cC$ closed under isomorphisms, say $\cE$ and $\cM$, such that morphisms in $\cE$ are left orthogonal to those in $\cM$, and each morphism $f$ in $\cC$ factors as $f=m\circ e$ where $e\in\cE$ and $m\in \cM$.

$Aut(n)$ is the set of automorphisms of $n$ in $\bT$. As in any Lawvere theory $\bT$, for $n\in \o$,  $n$ is canonically isomorphic to $1^n$ we always have a function
\[ \rho_n : S_n \times Aut(1)^n \lra Aut(n) \]
such that
\[ (\sigma, a_1,\ldots, a_n) \mapsto a_1\times  \ldots \times a_n \circ \pi_\sigma \]
i.e. $\rho_n$ sends a permutation $\sigma$ and $n$ isomorphisms of $1$ to an isomorphism of $n$ in $\bT$. We say that $\bT$ has {\em simple automorphisms} iff $\rho_n$ is a bijection, for $n\in\o$. Clearly, if $\bT$ has simple automorphisms then $2$ is not initial in $\bT$.

A Lawvere theory $\bT$ is {\em analytic} iff  structural morphisms and analytic morphisms form a factorization system in $\bT$ and  $\bT$ has simple automorphisms. A Lawvere theory $\bT$ is {\em rigid} iff it is analytic and the symmetric groups $S_n$ acting on $\bT(n,1)$ by permuting factors act freely on analytic morphisms, for $n\in \o$.

An {\em analytic interpretation} of Lawvere theories is an interpretation of Lawvere theories that preserves analytic morphisms. Thus we have a non-full subcategory of analytic Lawvere theories and analytic interpretations $\ALT$.  The latter has as a full subcategory the category $\RiLT$ of rigid Lawvere theories. We have inclusion functors
\[ \RiLT\lra \ALT \lra \LT \]
with the first one being a full inclusion.

We have an easy

\begin{lemma} \label{factorization}
In any analytic Lawvere theory $\bT$ any morphism $f:n\ra m$ has a factorization
\begin{center} \xext=800 \yext=450
\begin{picture}(\xext,\yext)(\xoff,\yoff)
  \settriparms[1`1`-1;400]
  \putVtriangle(0,0)[n`m`k;f`\pi_\phi`a]
\end{picture}
\end{center}
with $a$ being an analytic morphism in $\bT$ and $\phi:(k]\ra (n]$ a function. Such a factorization is unique up to a permutation, that is if $f= a'\circ \pi_{\phi'}$ is another such factorization there is $\sigma \in S_k$ such that
\[ \phi\circ \sigma = \phi',\;\;\;\; a = a'\circ \pi_\sigma.\]
\end{lemma}

{\em Proof.} When $\bT$ has simple automorphisms any structural morphism $s:n\ra m$ in $\bT$ can be presented as $(a_1\times\ldots, a_m)\circ \pi_\phi$ for some function $\phi :(m]\ra (n]$ and $a_i\in Aut(1)$ for $i\in (m]$. Thus if $f=a\circ s$ is a structural-analytic factorization of $f$, then $f=(a\circ(a_1\times\ldots, a_m))\circ \pi_\phi$ is also one. $\boxempty$

\subsection*{Monads} We shall consider three categories of finitary monads on $Set$. The category of all finitary monads with usual morphisms of monads will be denoted by $\Mnd$. A morphism of monads  $\tau: (M,\eta,\mu)\ra (M',\eta',\mu')$ is a natural transformation $\tau: M \ra M'$ such that $\tau\circ \eta^M=\eta^{M'}$ and $\tau\circ \mu^M=\mu^{M'}\circ \tau_{M'}\circ M(\tau)$.

Recall that a finitary monad $(M,\eta,\mu)$ on $Set$ is {\em analytic} iff $M$ weakly preserves wide pullbacks and both $\eta$ and $\mu$ are weakly cartesian natural transformations. A {\em morphism of analytic monads} on $Set$ $\tau: (M,\eta,\mu)\ra (M',\eta',\mu')$ is a weakly cartesian natural transformation $\tau$ that is a morphism of monads, (cf. \cite{J}, \cite{Z1}). Recall that a finitary monad $(M,\eta,\mu)$ is a {\em polynomial monad} on $Set$ iff $M$ preserves wide pullbacks and both $\eta$ and $\mu$ are cartesian natural transformations. Both types of functors and monads have much more explicit description (cf. \cite{J}, \cite{Z1}).

The categories of analytic and  polynomial monads with the suitable morphisms will be denoted by $\AMnd$ and $\PMnd$, respectively.
We have two inclusion functors
\[  \PMnd \lra \AMnd \lra \Mnd \]
the first one being full (cf. \cite{Z1}), the second full on isomorphisms.

\subsection*{Operads}
The symmetric operads provide yet another way of presenting models of equational theories. This kind of presentation is usually very convenient, but the models defined by such operads are more specific. The precise characterization of this more specific situation is the main objective of the paper.

Recall that a {\em symmetric operad} $\cO$ consists a family of sets $\cO_n$, for $n\in \o$, a unit element $\iota\in \cO_1$, for any $k,n,n_1,\ldots, n_k \in \o$ with $n={\sum_{i=1}^k n_i}$, a composition operation
\[  \ast : \cO_{n_1}\times\ldots \times \cO_{n_k} \times  \cO_k \lra \cO_n\]
a left action of the symmetric groups
\[  \cdot : S_n\times \cO_n \lra \cO_n \]
for $n\in \o$, such that the composition is associative with unit $\iota$ and compatible with group actions. A {\em morphism of symmetric operads} $f: \cO \ra \cO'$ is a function that respects arities of operations, unit, compositions, and group actions. For more on symmetric operads and their history one can consult for example \cite{Le}, but there the symmetric group act on the right.

The symmetric operad of symmetries $\Sym$ is defined as follows. The set of $n$-ary operations of $\Sym$ is the symmetric group $S_n$ on which $S_n$ act on the left by multiplication. The composition
\[  \star : S_{n_1}\times\ldots \times S_{n_k} \times  S_k \lra S_n\]
for $(\sigma_1,\ldots,\sigma_k;\tau) \in  S_{n_1}\times\ldots \times S_{n_k} \times  S_k$ the permutation
\[  \lk\sigma_1,\ldots,\sigma_k\rk\star\tau : n=\sum_{i=1}^k n_{\tau(i)} \lra n=\sum_{i=1}^k n_i  \]
is given by
\[ \lk i,r\rk \mapsto \lk \tau(i),\sigma_{\tau(i)}(r) \rk \]
where we consider the obvious lexicographic order on both $\sum_{i=1}^k n_{\tau(i)}$ and $\sum_{i=1}^k n_{\tau(i)}$. Note that even if composition is a functions between groups it is not a homomorphism of groups in general.

The category of rigid operads\footnote{The (colored) rigid operads where introduced in \cite{HMP} as `multicategories with non-standard amalgamation'. As this notion turned to be important in many other contexts it deserves a shorter name. The name we have chosen is related to the characterization we are going to prove in Theorem \ref{eq_rigid-the_polymonad}.} $\RiOp$ can be identified with the full subcategory of symmetric operads whose objects are those operads that have all the actions of symmetric groups free. For more the reader can consult \cite{HMP} and \cite{Z1}.

We denote by $\SOp$, $\RiOp$ the categories of symmetric and rigid operads, respectively. We have the `symmetrization' functor (see \cite{Z1})
\[ \P: \RiOp \lra \SOp \]
that we identify here with an embedding of a full subcategory.

%\newpage
\section{The equivalence of the three approaches}\label{3equivalences}
We shall recall the functors that exhibit equivalences of the following three categories $\ET$, $\LT$ and $\Mnd$:

\begin{center} \xext=1600 \yext=400
\begin{picture}(\xext,\yext)(\xoff,\yoff)
\putmorphism(800,200)(1,0)[\phantom{\LT}`\Mnd`\cM_l]{800}{1}a
\putmorphism(800,150)(1,0)[\phantom{\LT}`\phantom{\Mnd}`\cL_m]{800}{-1}b
  \putmorphism(0,200)(1,0)[\ET`\LT`\cL_e]{800}{1}a
   \putmorphism(0,150)(1,0)[\phantom{\ET}`\phantom{\Mnd}`\cE_l]{800}{-1}b
\end{picture}
\end{center}
As we described it in the introduction, the names of the functors are so chosen to remember their codomains with indices remembering their domains. We often drop the indices when it does not lead to a confusion.

Each of the above categories comes equipped with a semantic functor associating to objects of those categories their categories of models. As all monads in $\Mnd$ are defined on $Set$ only, we consider the models only in $Set$. It is well known that the equivalences that we describe below respect those semantic functors.

\subsection*{The functor $\cL_e=\cL : \ET \lra \LT$}

Let $T=(L,A)$ be an equational theory. A morphism $n\ra m$ in $\cL(T)$ is an $m$-tuple $\lk [t_1:\vec{x}^n], \ldots, [t_m:\vec{x}^n] \rk : n \ra m$ where $[t_i:\vec{x}^n]$ is an equivalence class of terms is context $\vec{x}^n$ modulo provable equivalence from axioms in $A$. The identity on $n$ is $\lk [x_1:\vec{x}^n], \ldots, [x_n:\vec{x}^n] \rk : n \ra n$. The composition is given by simultaneous substitution as follows
\begin{center} \xext=2000 \yext=400
\begin{picture}(\xext,\yext)(\xoff,\yoff)
\putmorphism(0,260)(1,0)[\phantom{n}`m`{\lk [t_i:\vec{x}^n]\rk}_{i\in (m]} ]{1000}{1}a
 \putmorphism(1000,260)(1,0)[\phantom{m}`\phantom{k}`{\lk [s_j:\vec{x}^m]\rk}_{j\in (k]} ]{1000}{1}a
 \putmorphism(0,100)(1,0)[\phantom{m}`\phantom{k}`{\lk [s_j(\lk x_i \backslash t_i \rk_{i\in (m]}):\vec{x}^m]\rk}_{j\in (k]} ]{2000}{1}b
 \putmorphism(0,170)(1,0)[n`k` ]{2000}{0}b
\end{picture}
\end{center}
The $i$-th projection on $1$ is $\pi^n_i = \lk [ x_i:\vec{x}^n] \rk$.

Let $I: T \ra T'$ be an interpretation. The functor $\cL(I)$ is defined on a morphism
\[ \lk [t_1:\vec{x}^n], \ldots, [t_m:\vec{x}^n] \rk : n \ra m \]
in $\cL(T)$ as
\[ \lk [\bar{I}(t_1):\vec{x}^n], \ldots, [\bar{I}(t_{m}):\vec{x}^n] \rk : n \lra m\]

A routine verification shows that $\cL$ is indeed a functor into $\LT$.

\subsection*{The functor $\cE_l=\cE : \LT \lra \ET$}

Let $\bT$ be a Lawvere theory. Then $\bT(n,1)$ is the set of  $n$-ary operations of the theory $\cE(\bT)$, for $n\in \o$. The set of axioms $\cE(\bT)$ contains a linear-regular axiom
\[ g(x_1, \ldots,x_n)=f(f_1(x_1,\ldots,x_{n_1}),\ldots,   f_k(x_{1+\sum_i^{k-1}n_i},\ldots,x_{n}))  : \vec{x}^n \]
for any morphisms $f, f_i, g$ in $\bT$ such that $f\circ( f_1\times\ldots\times f_k ) = g$ holds in $\bT$, and a linear axiom
\[ x_i=\pi^n_i (x_1, \ldots,x_n) : \vec{x}^n \]
for any $n\in \o$ and $i\in (n]$. An interpretation of Lawvere theories $F:\bT\ra \bT'$ induces an interpretation of equational theories $\cE(F)$ such that
\[ \cE(F)(f) = F(f)(x_1, \ldots, x_n): \vec{x}^n \]
for $f:n\ra 1$ in $\bT$.

 \subsection*{The functor $\cL_m=\cL : \Mnd \lra \LT$}

For a monad $M=(M,\eta, \mu)$, the category $\cL(M)$ is the dual of the full subcategory of the Kleisli category for $M$, spanned by the natural numbers.  In detail, for a monad $M$ we define the hom's in the category $\cL(M)$ as
\[ \cL(M)(n,m)= Set((m],M((n])) \]
for $n,m\in \o$. The compositions and identities are like in Kleisli category. The projection
\[ \pi^n_i : (1] \lra M((n])\]
sends $1$ to $\eta_{(n]}(i)$, for $n\in \o$ and $i\in (n]$.

For a morphism of monads $\tau : (M,\eta, \mu)\lra (M',\eta',\mu')$ and a morphism $f:n\ra m$ in $\cL(M)$ we put
\[ \cL(\tau)(f) (i) =\tau_{(n]}(f(i))\]
for $i\in (m]$.

\subsection*{The functor $\cM_l=\cM : \LT \lra \Mnd$}
For a Lawvere theory $\bT$, we define the monad $\cM(\bT)$ using coends. We put
\[ \cM(\bT) (X) = \int^{n\in \mathds{F}} X^n\times \bT(n,1) \]
for $X\in Set$.
The unit of $\cM(\bT)$ is
\[ \eta^\bT_X : X \ra \cM(\bT) (X) \]
sends $x\in X$ to the class of the element $\lk id_1,\bar{x} \rk$ where $id_1$ is the identity on $1$ in $\bT$ and $\bar{x} :(1]\ra X$ is the function picking $x$, i.e. $\bar{x}(1)=x$. The iterated functor $\cM^2(\bT)$ is given, for $X$ in $Set$ by
\[ \cM^2(\bT) (X)=\int^{m,n_1,\ldots,n_m \in \mathds{F}} X^{n}\times \bT(n_1,1)\times\ldots\times \bT(n_m,1)\times\bT(m,1) \]
where $n=\sum_{i=1}^m n_i$.
The multiplication of the monad $\cM(\bT)$
\[ \mu^\bT_X : \cM^2(\bT) (X) \lra \cM(\bT) (X) \]
is defined on components
\[  X^n\times \bT(n_1,1)\times\ldots \bT(n_m,1)\times \bT(m,1)   \lra  X^n\times \bT(n,1)\]
by composition, i.e. for $f:m\ra 1, f_1:n_1\ra 1,\ldots , f_m:n_m\ra 1$ in $\bT$ and $\vec{x}:(n]\ra X$
\[ \mu^\bT_X(\vec{x}, f_1,\ldots , f_m,f) = \lk \vec{x}, f\circ ( f_1\times\ldots \times f_m) \rk \]
where again $n=\sum_{i=1}^m n_i$.

%\newpage
\section{Lawvere theories vs Operads}

In this section we study the relations between Lawvere theories and operads, both symmetric and rigid. We shall
describe the adjunction  $\cP_a\dashv \cQ_f$ and the properties of the embeddings $\cP_a$ and $\cP_p$.
\begin{center} \xext=1600 \yext=450
\begin{picture}(\xext,\yext)(\xoff,\yoff)
   \putmorphism(800,100)(1,0)[\SOp`\LT`]{800}{0}a
   \putmorphism(800,50)(1,0)[\phantom{\SOp}`\phantom{\LT}`\cQ_f]{800}{-1}b
   \putmorphism(800,150)(1,0)[\phantom{\SOp}`\phantom{\LT}`\cP_a]{800}{1}a
   \putmorphism(0,100)(1,0)[\RiOp`\phantom{\SOp}`\cP]{800}{1}a
       \put(700,360){$\P_p$}
  \put(0,200){\line(0,1){120}} \put(0,320){\line(1,0){1600}} \put(1600,320){\vector(0,-1){120}}
\end{picture}
\end{center}

\subsection*{The functor $\cP_a : \SOp\ra \LT$}
Let $\cO$ be a symmetric operad: $\iota$, $\cdot$, $\ast$ denote the unit, symmetric groups actions, and compositions in $\cO$, respectively.
We define a Lawvere theory $\cP_a(\cO)$ as follows. The set of objects of $\cP_a(\cO)$ is the set of natural numbers $\o$. A morphism from $n$ to $m$ is an equivalence class of spans
\begin{center} \xext=800 \yext=520
\begin{picture}(\xext,\yext)(\xoff,\yoff)
 \settriparms[1`1`0;400]
  \putAtriangle(0,0)[r`n`m;\phi`\lk f,g_i\rk_{i\in m}`]
\end{picture}
\end{center}
such that $\phi:(r]\ra (n]$ is a function, $f:(r]\ra (m]$ is a monotone function, $r_i=|f^{-1}(i)|$ and we have $g_i\in \cO_{r_i}$ for $i\in (m]$ and $r=\sum_{i=1}^m r_i$. Two spans $\lk \phi,f,g_i\rk_{i\in m}$ and $\lk \phi ',f',g'_j\rk_{j\in m'}$ are equivalent iff $f=f'$, and there are permutations $\sigma_i\in S_{r_i}$ for $i\in (m]$
\begin{center} \xext=800 \yext=920
\begin{picture}(\xext,\yext)(\xoff,\yoff)
 \settriparms[1`1`0;400]
  \putAtriangle(0,400)[r`n`m;\phi`\lk f,g_i\rk_{i}`]
  \settriparms[0`-1`-1;400]
  \putVtriangle(0,0)[\phantom{n}`\phantom{m}`r';`\phi'`{\lk f',g'_i\rk_{i}}]

  \putmorphism(400,750)(0,-1)[``]{700}{-1}a
    \put(420,380){$_{\sum_i \sigma_i}$}
\end{picture}
\end{center}
such that
\[  g_i  = \sigma_i \cdot g'_i,\;\;\;\; \phi\circ \sum_i \sigma_i =\phi'\]
By $\sum_i \sigma_i:r \ra r$ we mean the permutation that is formed by placing permutations $\sigma_i$ `one after another'. Thus, it respects the fibers of $f$, i.e.  $f\circ \sum_i \sigma_i =f$. Clearly, we shall deal with the spans when we perform constructions on morphisms in $\cP_a(\cO)$, but when we consider equalities between spans we shall invoke the above equivalence relation.

The composition $\lk \phi'',f'',g''_j\rk_{i\in (k]}: n \ra k$ of two morphism $\lk \phi,f,g_i\rk_{i\in (m]}: n \ra m$ and $\lk \phi ',f',g'_j\rk_{j\in (k]}:m\ra k$ is defined as follows. In the diagram
\begin{equation}\label{compo}\end{equation}
\begin{center} \xext=2000 \yext=720
\begin{picture}(\xext,\yext)(\xoff,\yoff)
 \settriparms[1`1`0;400]
  \putAtriangle(200,0)[r`n`m;\phi``]
   \putAtriangle(1000,0)[r'`\phantom{m}`k;\phi'``]
   \put(580,180){$_{\lk f,g_i\rk_{i}}$}
    \putAtriangle(600,400)[r''`\phantom{r}`\phantom{r'};\bar{\phi}`\bar{f}`]
   \put(1350,180){$_{\lk f',g'_j\rk_{j}}$}

   \put(800,900){\line(-1,-1){750}}
   \put(800,900){\line(1,-1){100}}
   \put(50,150){\vector(1,-1){100}}
   \put(250,480){$\phi''$}

    \put(1200,900){\line(1,-1){750}}
   \put(1200,900){\line(-1,-1){100}}
   \put(1950,150){\vector(-1,-1){100}}
   \put(1650,480){$\lk f'',g''_j \rk_j$}
\end{picture}
\end{center}
the square is a pullback of $f$ along $\phi'$. The function $\bar{f}$ is chosen so that it is monotone. We put $f''=f'\circ \bar{f}$,  $\phi''=\phi\circ \bar{\phi}$, and $g_j''= g_j'\ast \lk g_{\phi(l)} \rk_{l\in f^{-1}(j)}$.

The identity on $n$ is the span
\begin{center} \xext=800 \yext=500
\begin{picture}(\xext,\yext)(\xoff,\yoff)
 \settriparms[1`1`0;400]
  \putAtriangle(0,0)[n`n`n;id_n`\lk id_n,\iota\rk_{i}`]
\end{picture}
\end{center}
As $S_1$ contains the identity permutation only, any span equivalent to an identity span is actually equal to it.

The projection $\pi^n_i:n\ra 1$ on $i$-th coordinate is the span
\begin{center} \xext=800 \yext=520
\begin{picture}(\xext,\yext)(\xoff,\yoff)
 \settriparms[1`1`0;400]
  \putAtriangle(0,0)[1`n`1;\bar{i}`\lk id,\iota\rk`]
\end{picture}
\end{center}
where $i\in (n]$ and $\bar{i}(1)=i$.

For a morphism of symmetric operads $h:\cO\ra \cO'$ we define a functor
\[ \cP_a(h): \cP_a(\cO) \lra \cP_a(\cO') \]
so that for a morphism $\lk \phi,f,g_i\rk_{i\in (m]}: n \ra m$ in $\cP_a(\cO)$ we define a morphism
\[ \cP_a(h)(\lk \phi,f,g_i\rk_{i\in (m]}) = \lk \phi,f,h(g_i)\rk_{i\in (m]}: n \ra m \]
in  $\cP_a(\cO')$.

This ends the definition of the functor $\P_a$.
\subsection*{The functor $\cQ_f : \LT \lra \SOp$}

Let $\bT$ be a Lawvere theory. The operad $\cQ_f(\bT)$ consists of operations of $\bT$, i.e. morphisms to $1$. In detail it can be described as follows. The set of $n$-ary operations $\cQ_f(\bT)_n$ is the set of  $n$-ary operations  $\bT(n,1)$ of $\bT$, for $n\in \o$. The action
\[ \cdot : S_n \times \cQ_f(\bT)_n \lra \cQ_f(\bT)_n \]
is given, for $f\in\bT(n,1)$ and $\sigma\in S_n$,  by
\[ \sigma \cdot  f = f\circ \pi_{\sigma} \]
The identity of $\cQ_f(\bT)$ is $\iota=id_1\in\bT(1,1)$. The composition
\[ \ast : \cQ_f(\bT)_{n_1}\times\ldots\times\cQ_f(\bT)_{n_{k}}\times\cQ_f(\bT)_k\lra \cQ_f(\bT)_n \]
is given, for $f\in  \cQ_f(\bT)_k$ and $f_i\in \cQ_f(\bT)_{n_i}$,  where $i\in (k]$, $n={\sum_{i\in k} n_i}$, by
\[ \lk f_1,\ldots , f_{k}\rk \ast f = f\circ (f_1\times\ldots ,\times f_{k}) \]
where $f_1\times\ldots ,\times f_{k}$ is defined using the chosen projections in $\bT$ and $\circ$ is the composition in $\bT$.

If $F:\bT\ra \bT'$ is a morphism of Lawvere theories then the map of symmetric operads
\[ \cQ_f(F) : \cQ_f(\bT) \ra \cQ_f(\bT') \]
is defined, for $f\in \cQ_f(\bT)_n$, by
\[ \cQ_f(F)(f) = F(f) \]

This ends the definition of the functor $\cQ_f$.

\subsection*{The adjunction $\cP_a\dashv \cQ_f$ and the properties of the functor $\cP_a$}

We note for the record
\begin{proposition} \label{funtors_P_Op}
The functors $\cP_a : \SOp\lra \LT$ and $\cQ_f : \LT \ra \SOp$ are well defined. $\boxempty$
\end{proposition}

We have an easy

\begin{lemma} \label{iso}
Let $\cO$ be a symmetric operad and $n\in\o$. An automorphism on $n$ in  $P_a(\cO)$ is represented by a span of the following form
\begin{center} \xext=800 \yext=500
\begin{picture}(\xext,\yext)(\xoff,\yoff)
 \settriparms[1`1`0;400]
  \putAtriangle(0,0)[n`n`n;\phi`\lk id_n,a_i\rk_{i}`]
\end{picture}
\end{center}
where $\phi:(n]\ra (n]$ is a bijection, $a_i\in \cO_1$ is an invertible operation, i.e. there is $b_i\in \cO_1$ such that $a_i\ast b_i=\iota =b_i\ast a_i$ for $i\in (n]$. It is the unique span in its equivalence class.
\end{lemma}

{\em Proof.}
Consider a pair of morphisms in $\P_a(\cO)$
\begin{center} \xext=2000 \yext=520
\begin{picture}(\xext,\yext)(\xoff,\yoff)
 \settriparms[1`1`0;500]
  \putAtriangle(0,0)[r`n`n;\phi``]
  \putAtriangle(1000,0)[r`\phantom{n}`n;\phi'`\lk f',h_i\rk_{i}`]
  \put(430,200){$\lk f,g_j\rk_{j}$}
\end{picture}
\end{center}
that are inverse  one to the other. As the above composition is an identity it follows that $\phi$ and $f'$ are epi. Thus, because of the other composition
$\phi'$ and $f$ are surjections, as well. As pulling back along a surjection reflects injections, all functions $\phi$, $f$, $\phi'$ and $f'$ must be also injective and hence bijective. Then it is easy to see that $g_{\phi'(j)}$ is an inverse of $h_j$ for $j\in (n]$.
$\boxempty$

\begin{proposition} \label{operad-emb}
We have an adjunction $\cP_a\dashv \cQ_f$. The functor $\cP_a$ is faithful.
\end{proposition}

{\em Proof.}
 We first show that $\cP_a\dashv \cQ_f$. For a symmetric operad $\cO$ the unit is
\[ \eta_\cO :\cO \lra \cQ_f (\cP_a(\cO)) \]
\[ \cO_n \ni g \mapsto \lk id_n,!,g\rk  \]
For Lawvere theory $\bT$ the counit is
\[ \varepsilon_\bT: \cP_a\cQ_f(\bT)\lra \bT   \]
\[ \lk \phi,f, g_i\rk_{i\in (m]} \mapsto (g_1\times \ldots \times g_{m})\circ \pi_\phi  \]

We verify the triangular equalities. For $g\in \cQ_f(\bT)_n=\bT(n,1)$ we have
\[ \cQ_f (\varepsilon_\bT) \circ \eta_{\cQ_f(\bT)}(g) = \]
\[  =\cQ_f (\varepsilon_\bT) (\lk id_n,!,g\rk) = \]
\[ = g\circ \pi_{id_n} = g\]
For $\lk \phi,f,g_i\rk_{i\in (m]}\in \cP_a(\cO)$ we have
\[ \varepsilon_{\cP_a(\cO)} \circ \cP_a(\eta_\cO)(\lk \phi,f,g_i \rk_{i\in (m]} ) = \]
\[ =\varepsilon_{\cP_a(\cO)} (\lk \phi,f,\lk id_{r_i},!,g_i\rk \rk_{i\in (m]} ) = \]
\[ = (\lk id_{r_1},!,g_1 \rk\times \ldots \times\lk id_{r_{m}},!,g_{m} \rk) \circ \pi_\phi  =\]
\[ = \lk \phi,f,g_i \rk_{i\in (m]}  \]
As the unit $\eta_\cO$ is mono, $\cP_a$ is faithful.
$\boxempty$

\begin{proposition}\label{operad-emb-im}
The functor $\cP_a$ is faithful, full on isomorphisms and its essential image is the category of analytic Lawvere theories $\ALT$ i.e. it factorizes as an equivalence of categories $\cL_o$ followed by $\cP_a^l$
\begin{center}
\xext=600 \yext=650
\begin{picture}(\xext,\yext)(\xoff,\yoff)
 \settriparms[-1`-1`-1;600]
 \putbtriangle(0,0)[\LT`\ALT`\SOp;\P_a^l`\cP_a`\cL_o]
\end{picture}
\end{center}
\end{proposition}

{\em Proof.} Recall that we have a unique morphism of Lawvere theories from the initial theory $\pi : \mathds{F}^{op}\ra\cP_a(\cO)$. For a function $\phi : (m]\ra (n]$, $\pi_\phi$ the morphism $\pi_\phi$ is represented by the span of the form
\begin{center} \xext=800 \yext=520
\begin{picture}(\xext,\yext)(\xoff,\yoff)
 \settriparms[1`1`0;400]
  \putAtriangle(0,0)[m`n`m;\phi`\lk id_m,\iota\rk_{i\in (m]}`]
\end{picture}
\end{center}

The class of the structural morphisms in $\cP_a(\cO)$ is the closure under isomorphism of the class of morphisms $\{\pi_\phi : \phi\in \mathds{F}  \}$. It is easy to see that the structural morphisms in $\cP_a(\cO)$ are (represented by) the spans of the form
\begin{center} \xext=800 \yext=520
\begin{picture}(\xext,\yext)(\xoff,\yoff)
 \settriparms[1`1`0;400]
  \putAtriangle(0,0)[m`n`m;\phi`\lk id_m,a_i\rk_{i\in (m]}`]
\end{picture}
\end{center}
where $\phi$ is any function and $a_i$ is an invertible unary operation, for $i\in (m]$. Thus by Lemma \ref{iso}, a morphism is an isomorphism in $\cP_a(\cO)$ iff it is represented by a span as above with $\phi$ being a bijection.

The analytic morphisms in $\cP_a(\cO)$ are (represented by) the
spans of the form
\begin{center} \xext=800 \yext=520
\begin{picture}(\xext,\yext)(\xoff,\yoff)
 \settriparms[1`1`0;400]
  \putAtriangle(0,0)[n`n`m;\phi`\lk f,g_i\rk_{i\in (m]}`]
\end{picture}
\end{center}
where $\phi$ is a bijections.

Clearly, both classes contain isomorphisms and are closed under composition.

Any morphism $\lk \phi,f,g_i\rk_{i\in (m]}:n\ra m$ in $\cP_a(\cO)$ has a structural-analytic factorization as follows
\begin{center} \xext=2000 \yext=520
\begin{picture}(\xext,\yext)(\xoff,\yoff)
 \settriparms[1`1`0;500]
  \putAtriangle(0,0)[r`n`r;\phi``]
  \putAtriangle(1000,0)[r`\phantom{r}`n;id_r`\lk f,g_i\rk_{i}`]
  \put(430,200){$\lk id_r,\iota\rk_{j}$}
\end{picture}
\end{center}

Thus to show that structural and analytic morphisms form a factorization system it remains to show that structural morphisms are left orthogonal to the analytic morphisms. Let
\begin{center} \xext=1250 \yext=1250
\begin{picture}(\xext,\yext)(\xoff,\yoff)
%lt
\setsqparms[-1`-1`0`0;600`500]
 \putsquare(50,600)[n`r`\phantom{k}`;\psi`\phi``]
%rt
\setsqparms[1`0`-1`0;600`500]
 \putsquare(650,600)[\phantom{r}`m``m;\lk f,h_i \rk_{i\in (m]}``1_m`]
%lb
\setsqparms[0`1`0`-1;600`500]
 \putsquare(50,100)[k``k`r';`\lk 1_k,a_j \rk_{j\in (k]}``\phi']
%rb
\setsqparms[0`0`1`1;600`500]
\putsquare(650,100)[`\phantom{m}`\phantom{r'}`1;``\lk !,g \rk`\lk !,g' \rk]

\putmorphism(650,1050)(0,-1)[``\sigma]{900}{1}r
\end{picture}
\end{center}
be a commutative square in $\cP_a(\cO)$ with left vertical morphism $\lk \phi,1_r,a_i\rk_{j\in (k]}$ being a structural map and right vertical morphism  $\lk 1_m,!,g\rk$ an analytic map. We have chosen the right bottom to be $1$ to simplify notation but the general case is similar. The commutation means that $r=r'$ and there is a permutation $\sigma\in S_r$  such that $$\psi=\phi\circ \phi'\circ \sigma$$ and
$$\lk a_{\phi'(1)},\ldots,a_{\phi'(r)} \rk \ast g' = \sigma\cdot  ( \lk h_1,\ldots,h_{m} \rk \ast g)$$
Putting into the square a diagonal morphism $\lk \phi'\circ \sigma,f,\bar{h_i}\rk_{i\in (m]}$
\begin{center} \xext=1250 \yext=1250
\begin{picture}(\xext,\yext)(\xoff,\yoff)
%lt
\setsqparms[-1`-1`1`0;600`500]
 \putsquare(50,600)[n`r`\phantom{k}`r;\psi`\phi`1_r`]
%rt
\setsqparms[1`0`-1`0;600`500]
 \putsquare(650,600)[\phantom{r}`m``m;\lk f,h_i \rk_{i}``1_m`]
%lb
\setsqparms[0`1`1`-1;600`500]
 \putsquare(50,100)[k`\phantom{r}`k`r;`\lk 1_k,a_j \rk_{j}`\sigma`\phi']
%rb
\setsqparms[0`0`1`1;600`500]
 \putsquare(650,100)[`\phantom{m}`\phantom{r}`1;``\lk !,g \rk`\lk !,g' \rk]

  \put(200,400){$_{\phi'\circ \sigma}$}
  \put(850,720){$_{\lk f,\bar{h_i} \rk_{i}}$}
  \put(600,550){\vector(-4,-3){500}}
   \put(700,650){\vector(4,3){500}}
\end{picture}
\end{center}
where
\[ \bar{h_i} = \lk a^{-1}_{\phi'\circ \sigma(l)} \rk_{l\in f^{-1}(i)}  \ast h_i \]
we see that the permutations $1_r$ and $\sigma$ show that both triangles commute. It is not difficult to see that this diagonal filling is unique.
Thus analytic morphisms are indeed right orthogonal to the structural ones and  $\cP_a(\cO)$ is an analytic Lawvere theory.

From the description of the functor  $\cP_a(h): \cP_a(\cO)\ra \cP_a(\cO')$ and the description of the structure of $\cP_a(\cO)$ it is clear that $\cP_a(h)$ sends the analytic (structural) morphisms to the analytic (structural) ones. Thus $\cP_a(h)$ is an analytic interpretation of Lawvere theories.

Now let $\bT$ be any Lawvere theory. As the class of analytic morphisms in $\bT$ is right orthogonal to a class of morphisms, it is closed under finite products and isomorphisms. In particular, a composition of an analytic morphism $f:n\ra 1$ in $\bT$ with a permutation morphism $\pi_\sigma$ with $\sigma\in S_n$ is again an analytic morphism. Thus the analytic operations of any Lawvere theory $\bT$ form a symmetric operad. The composition $\lk f_1, \ldots, f_n\rk \ast f$ is defined to be $f\circ (f_1\times \ldots \times f_n)$ and the action of $\sigma\in S_n$ on an analytic morphism $f:n\ra 1$ is $\sigma \cdot f = f\circ \pi_\sigma$. The unit is the identity morphism on $1$. So defined the symmetric part of the operad $\bT$ will be denoted as $\bT^s$. We have an inclusion morphism of symmetric operads $$\bT^s\ra \cQ_f(\bT)$$ By adjunction we get a morphism
\[ \psi_\bT:\cP_a(\bT^s) \lra \bT \]
Clearly, $\psi_\bT$ is bijective on objects. If $\bT$ is analytic then $\psi_\bT$ is full (faithful) since the structural-analytic factorization exists (is unique and $\pi:\F\ra \bT$ is faithful), see Lemma \ref{factorization}.

If $I:\bT\ra \bT'$ is an analytic interpretation between any Lawvere theories, then the diagram
\begin{center} \xext=1000 \yext=600
\begin{picture}(\xext,\yext)(\xoff,\yoff)
\setsqparms[1`1`1`1;1000`500]
 \putsquare(0,50)[\P(\bT^s)`\P(\bT'^s)`\bT`\bT;\P(I^s)`\psi_\bT`\psi_{\bT'}`I]
 \end{picture}
\end{center}
commutes, where $I^s$ is the obvious restriction of $I$ to $\bT^s$. Thus the essential image of $\cP_a$ is indeed the category of analytic Lawvere theories and analytic interpretations. An isomorphic interpretation of Lawvere theories is always analytic. Therefore $\cP_a$ is full on isomorphisms. $\boxempty$

We have
\begin{proposition}\label{monadic}
The functor $\cQ_f : \LT \ra \SOp$ is monadic.
\end{proposition}
{\em Proof.}  We shall verify that $\cQ_f$ satisfies the assumptions of Beck monadicity theorem. By Proposition \ref{operad-emb}, $\cQ_f$ has a left adjoint. It is easy to see that $\cQ_f$ reflects isomorphisms. We shall verify that $\LT$ has and $\cQ_f$ preserves $\cQ_f$-contractible  coequalizers.

Let $I,I': \bT' \ra \bT$ be a pair of interpretations between Lawvere theories so that
\begin{center} \xext=1400 \yext=400
\begin{picture}(\xext,\yext)(\xoff,\yoff)
\putmorphism(800,200)(1,0)[\phantom{\cQ_f(\bT)}`\cO`q]{600}{1}a
\putmorphism(800,150)(1,0)[\phantom{\cQ_f(\bT)}`\phantom{\cO}`s]{600}{-1}b
  \putmorphism(0,200)(1,0)[\cQ_f(\bT')`\cQ_f(\bT)`\cQ_f(I')]{800}{1}b
   \putmorphism(0,50)(1,0)[\phantom{\cQ_f(\bT')}`\phantom{\cQ_f(\bT)}`r]{800}{-1}b
    \putmorphism(0,280)(1,0)[\phantom{\cQ_f(\bT')}`\phantom{\cQ_f(\bT)}`\cQ_f(I)]{800}{1}a
\end{picture}
\end{center}
is a split coequalizer in $\SOp$.  We define a Lawvere theory $\bT_\cO$ so that a morphism from $n$ to $m$ in $\bT_\cO$ is an $m$-tuple $\lk g_1, \ldots, g_m\rk$ with $g_i\in \cO_n$, for $i=1,\ldots, m$. The compositions and the identities in $\bT_\cO$ are defined in the obvious way from the compositions and the unit in $\cO$. The projections $\bar{\pi}^n_i$ in $\bT_\cO$ are the images of the projections $\pi^n_i$ in $\bT$, i.e. $\bar{\pi}^n_i= q(\pi^n_i)$.

The functor $\tilde{q}: \bT \ra \bT_\cO$ is defined, for $f:n\ra m$ in $\bT$, as
\[ \tilde{q}(f) = \lk q(\pi^m_1\circ f),\ldots, q(\pi^m_m\circ f)\rk \]

First we verify, that $\bT_\cO$ has finite products. For this, it is enough to verify that $ \lk f_1, \ldots, f_n\rk\ast\bar{\pi}^n_i  = f_i$, where $\ast$ is the composition in the operad $\cO$. The uniqueness of the morphism into the product is obvious from the construction. We have routine calculations
\[  \lk f_1, \ldots, f_n\rk \ast\bar{\pi}^n_i  = \]
\[  q\circ s( \lk f_1, \ldots, f_n\rk \ast q(\pi^n_i) ) = \]
\[  \lk q\circ s(f_1), \ldots, q \circ s(f_n)\rk\ast (q\circ s\circ q(\pi^n_i)) = \]
\[  \lk q\circ s(f_1), \ldots,q\circ s(f_n)\rk \ast (q\circ Op(I)\circ r(\pi^n_i)) = \]
\[  \lk q\circ s(f_1), \ldots,q\circ s(f_n)\rk  \ast (q\circ Op(I')\circ r(\pi^n_i))  = \]
\[   \lk q\circ s(f_1), \ldots,q\circ s(f_n)\rk \ast (q(\pi^n_i)) = \]
\[  q(\lk  s(f_1), \ldots,s(f_n)\rk\ast\pi^n_i) = \]
\[  q(s(f_i)) = f_i \]

It is obvious that $\tilde{q}$ is a morphism of Lawvere theories and that $\cQ_f(\tilde{q})=q$. It remains to verify that $\tilde{q}$ is a coequalizer in $\LT$.
Let $p: \bT \ra \bS$ be a morphism in $\LT$ coequalizing $I$ and $I'$.
\begin{center} \xext=1100 \yext=700
\begin{picture}(\xext,\yext)(\xoff,\yoff)
 \settriparms[1`1`1;500]
  \putqtriangle(600,0)[\bT`\bT_\cO`\bS;\tilde{q}`p`\tilde{k}]

  \putmorphism(0,500)(1,0)[\bT'`\phantom{\bT}`]{600}{0}b
   \putmorphism(0,450)(1,0)[\phantom{\bT'}`\phantom{\bT}`I']{600}{1}b
    \putmorphism(0,550)(1,0)[\phantom{\bT'}`\phantom{\bT}`I]{600}{1}a
\end{picture}
\end{center}
The morphism $\cQ_f(p)$ is coequalizing  $\cQ_f(I)$ and $\cQ_f(I')$ in $\SOp$. Thus there is a unique morphism $k$ in $\SOp$ making the triangel on the right
\begin{center} \xext=1300 \yext=700
\begin{picture}(\xext,\yext)(\xoff,\yoff)
 \settriparms[1`1`1;500]
  \putqtriangle(800,0)[\cQ_f(\bT')`\cO`\cQ_f(\bS);q`\cQ_f(p)`k]

  \putmorphism(0,500)(1,0)[\cQ_f(\bT')`\phantom{\cQ_f(\bT)}`]{800}{0}b
   \putmorphism(0,450)(1,0)[\phantom{\cQ_f(\bT')}`\phantom{\cQ_f(\bT)}`\cQ_f(I')]{800}{1}b
    \putmorphism(0,550)(1,0)[\phantom{\cQ_f(\bT')}`\phantom{\cQ_f(\bT)}`\cQ_f(I)]{800}{1}a
\end{picture}
\end{center}
commute. We define the functor $\tilde{k}$ so that
\[  \tilde{k}(\lk f_1,\ldots, f_n\rk) = \lk k(f_1),\ldots, k(f_n)\rk \]
for any morphism $\lk f_1,\ldots, f_n\rk$ in $\bT_\cO$. The verification that  $\tilde{k}$ is the required unique functor is left for the reader.  $\boxempty$

\subsection*{The functor $\cP^{ol}_{p}=\cP_{p} : \RiOp \ra \LT$}

The functor $\cP_p$ is defined as the composition of the functors $\cP_a\circ \cP$. We have

\begin{proposition}\label{im_rigid_in_lt}
The essential image of the functor $\cP_{p} : \RiOp\lra \LT$ is the category of $\RiLT$ of rigid Lawvere theories and analytic morphisms between them.
\end{proposition}
{\em Proof.} As $\cP$ is full and faithful, $\cP_p$ is faithful and full on analytic morphisms. The image of $\cP$ consists of those symmetric operads for which the symmetric group actions are free. Thus the image of $\cP_p$ consists of those analytic Lawvere theories in which the symmetric actions are free on analytic operations, i.e. it consists of the rigid Lawvere theories. $\boxempty$

 We end this section pointing out to yet another property of analytic Lawvere theories. Let $\bT$ be a category with finite products. A morphism $p:n\ra m$ in $\bT$ is a projection iff there is a morphism $p':n\ra m'$ so that the diagram
\begin{center} \xext=800 \yext=200
\begin{picture}(\xext,\yext)(\xoff,\yoff)
 \putmorphism(0,0)(1,0)[m`n`p]{400}{-1}a
 \putmorphism(400,0)(1,0)[\phantom{n}`m'`p']{400}{1}a
\end{picture}
\end{center}
is a product in $\cT$. We call such a diagram a {\em decomposition} of $n$. A decomposition is {\em trivial} iff $m$ or $m'$ is the terminal object (i.e. $0$ if $\bT$ is a Lawvere theory), otherwise it is non-trivial. An object is {\em indecomposable} if it does not have a non-trivial decomposition.

\begin{proposition}
$1$ is indecomposable in any analytic Lawvere theory.
\end{proposition}

{\em Proof.} It is enough to show that for any symmetric operad $\cO$, $1$ is indecomposable in $\cP_a(\cO)$. Consider the following diagram
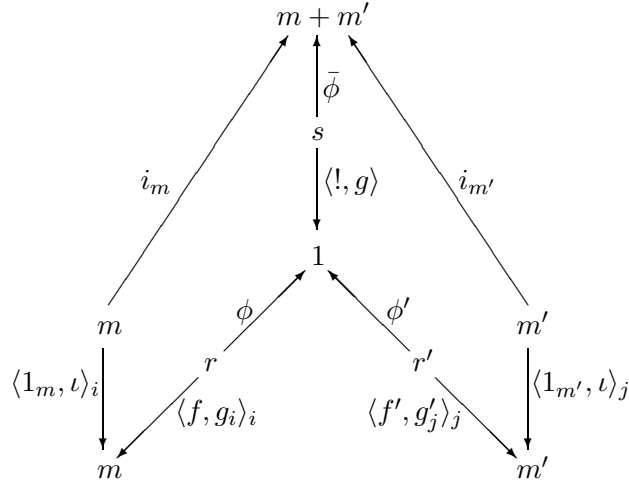
\begin{center} \xext=1600 \yext=1620
\begin{picture}(\xext,\yext)(\xoff,\yoff)
% 2 legs down
  \put(50,0){$m$}
  \put(1500,0){$m'$}
  \put(420,370){$r$}
  \put(1140,370){$r'$}
  \put(790,740){$1$}
  \put(400,350){\vector(-1,-1){280}}
  \put(320,200){$\lk f,g_i \rk_i$}
  \put(490,440){\vector(1,1){280}}
  \put(530,550){$\phi$}
  \put(1130,440){\vector(-1,1){280}}
  \put(1050,550){$\phi'$}
  \put(1220,350){\vector(1,-1){280}}
  \put(980,200){$\lk f',g'_j \rk_j$}
%vertical
    \put(810,1140){\vector(0,-1){280}}
    \put(830,1000){$\lk !,g\rk$}
  \put(810,1250){\vector(0,1){280}}
  \put(830,1320){$\bar{\phi}$}
    \put(790,1170){$s$}
    \put(670,1570){$m+m'$}
% sides
     \put(70,450){\vector(0,-1){350}}
     \put(-250,300){$\lk 1_m,\iota\rk_i$}
  \put(90,600){\vector(2,3){620}}
  \put(200,1000){$i_m$}
  \put(50,500){$m$}

   \put(1540,450){\vector(0,-1){350}}
   \put(1300,1000){$i_{m'}$}
  \put(1540,600){\vector(-2,3){620}}
  \put(1550,300){$\lk 1_{m'},\iota\rk_j$}
  \put(1500,500){$m'$}
 \end{picture}
\end{center}
We assume that the morphisms $\lk \phi,f,g_i\rk_i$, $\lk \phi',f',g'_j\rk_j$ are projections making $1$ into a product in $\cP_a(\cO)$. We also have two canonical projections from $m+m'$ to $m$ and $m'$. The morphism $\lk \bar{\phi},!,g\rk$ is a morphism into the product making both triangle commute.

From the commutations of the triangles easily follows that $$g_i\ast\lk g,\ldots,g\rk =\iota=g'_j\ast\lk g,\ldots,g\rk$$ for $i\in (m]$ and $j\in (m']$. This means that $g_i=g_j'=g^{-1}\in \cO_1$ for $i\in (m]$ and $j\in (m']$ and hence $r=m$, $r'=m'$, $f=1_m$, $f'=1_{m'}$. Moreover, $s=1$ and $!=1_1$.
Now commutativity says that there are $\sigma\in S_m$ and $\sigma'\in S_{m'}$ such that $i_m\circ \sigma = \bar{\phi}\circ \phi$ and $i_{m'}\circ \sigma' = \bar{\phi}\circ \phi'$. This is possible only if $m+m'=1$. $\boxempty$

From this Proposition follows immediately that the Lawvere theory of Jonsson-Tarski algebras is not analytic.

%\newpage
\section{Finitary Monads vs  Operads}

First we explain the diagram
\begin{center} \xext=800 \yext=1100
\begin{picture}(\xext,\yext)(\xoff,\yoff)
\setsqparms[1`-1`-1`0;800`500]
\putsquare(0,550)[\LT`\Mnd`\phantom{\SOp}`\phantom{\AMnd};\cM_l`\cP_a=\cP_a^{ol}`\cP_a^m`]
\setsqparms[1`-1`-1`1;800`500]
\putsquare(0,50)[\SOp`\AMnd`\RiOp`\PMnd;\cM_o^a`\cP=\cP_p^o`\cP_p^m`\cM_o^p]
 \end{picture}
\end{center}
commuting up to isomorphisms, with $\cP_a^m$ and $\cP_p^m$ being inclusions and $\cM_l$ is the equivalence of categories defined in \ref{3equivalences}. The remaining two horizontal functors are also equivalences of categories. We recall them below (cf. \cite{Z1}).

For a set $X$ we consider $X^n$ as the set of functions $X^{(n]}$. Then the permutation group $S_n$ acts naturally of $X^n$ on the right by composition. For a symmetric operad $\cO$, the monad $\cM_o^a(\cO)$ on a set $X$ is defined as
\[ \cM_o^a(\cO)(X)=\sum_{n\in \o} X^n\otimes_{n} \cO_n \]
Thus in $X^n\otimes_{n} \cO_n$ we identify $\lk\vec{x}\circ \sigma,f\rk$ with  $\lk \vec{x}, \sigma\cdot f\rk$ for $f\in \cO_n$, $\vec{x}: (n]\ra X$ and $\sigma\in S_n$.

For a rigid operad $\cO$ the monad  $\cM_o^p(\cO)$ on a set $X$ is defined as
\[ \cM_o^p(\cO)(X)=\sum_{n\in \o} X^n\times \cO_n \]
For more detailed description see for example \cite{Z1}. One can also find there the commutation of the lower square in the above diagram.

The commutation of the upper square is the content of the following.

\begin{proposition}
The square of categories and functors
\begin{center} \xext=800 \yext=600
\begin{picture}(\xext,\yext)(\xoff,\yoff)
\setsqparms[1`-1`-1`1;800`500]
\putsquare(0,50)[\LT`\Mnd`\SOp`\AMnd;\cM_l`\cP_a`\cP_a^m`\cM_o^a]
 \end{picture}
\end{center}
commutes up to an isomorphism.
\end{proposition}

{\em Proof.} Let $\cO$ be a symmetric operad. We need to define a natural isomorphism $\kappa$, so that
\[ \kappa^\cO : \cM_o^a(\cO) \lra \cM_l\cP_a(\cO) \]
is an isomorphism of monads natural in $\cO$.  The component of $\kappa^\cO$ at a set $X$
\[ \kappa^\cO_X : \sum_{n\in \o} X^n\otimes \cO_n \lra \int^{n\in \F} X^n\times \cP(\cO)(n,1) \]
is given by
\[  [\vec{x},a] \mapsto [\vec{x}, (1_n,!,a) ]   \]
where $\vec{x}:(n]\ra X$, $a\in \cO_n$ and $(1_n,!,a)$ is a span
\begin{center} \xext=800 \yext=520
\begin{picture}(\xext,\yext)(\xoff,\yoff)
 \settriparms[1`1`0;400]
  \putAtriangle(0,0)[n`n`1;1_n`\lk !,a\rk`]
\end{picture}
\end{center}
The verification that so defined $\kappa$ is indeed a natural isomorphism is left for the reader.
  $\boxempty$

\subsection*{The functor $\cQ_f^m : \Mnd \ra \AMnd$}

As the horizontal functors in the above diagram are equivalences of categories it follows from Proposition \ref{monadic} that the embedding functor $i:\AMnd \ra \Mnd$ has a right adjoint
\[ \cQ_f^m : \Mnd \ra \AMnd \]
which is monadic. In other words, any finitary monad on $Set$ is an  algebra for a monad on the category of analytic monads. We could define the functor $\cQ_f^m$ and the related monad $\Vb$ on $\AMnd$ directly, but we shall derive it from the more fundamental situation.

Let $\beta:\B\ra \F$ be the inclusion functor. It induces the following diagram of categories and functors that we describe below
\begin{center} \xext=3100 \yext=1300
\begin{picture}(\xext,\yext)(\xoff,\yoff)
%squares
\setsqparms[1`0`0`1;1600`750]
 \putsquare(150,450)[\hskip 25mm\Mnd\; =\; \Mon(\End)`\End`\hskip 25mm\AMnd\; =\; \Mon(\An)`\An;\widehat{U}```U]
\setsqparms[-1`0`0`-1;800`750]
 \putsquare(1750,450)[\phantom{\End}`Set^\F`\phantom{\An}`Set^\B;i_\F```i_\B]

%vertical
 \putmorphism(700,1150)(0,-1)[\phantom{\End}`\phantom{\An}`\Mon((-)^a)]{700}{1}r
  \putmorphism(600,1150)(0,-1)[\phantom{\End}`\phantom{\An}`]{700}{-1}l

 \putmorphism(200,1150)(0,-1)[\phantom{\Mnd}`\phantom{2\bCat_{\times}}`\cQ_f^m]{700}{1}r
 \putmorphism(100,1150)(0,-1)[\phantom{\Mnd}`\phantom{2\bCat_{\times}}`\cP_a^m]{700}{-1}l

 \putmorphism(1800,1150)(0,-1)[\phantom{\End}`\phantom{\An}`(-)^a]{700}{1}r
 \putmorphism(1700,1150)(0,-1)[\phantom{\End}`\phantom{\An}`i^a]{700}{-1}l

  \putmorphism(2600,1150)(0,-1)[\phantom{Set^\F}`\phantom{Set^\B}`\beta^*]{700}{1}r
 \putmorphism(2500,1150)(0,-1)[\phantom{Set^\F}`\phantom{Set^\B}`Lan_\beta]{700}{-1}l

 \putmorphism(3100,1150)(0,-1)[\F`\B`\beta]{700}{-1}r

%\Vb=\Mon(\V)
  \put(150,180){\oval(100,100)[b]}
  \put(100,180){\line(0,1){200}}
  \put(200,180){\vector(0,1){200}}
  \put(0,0){${(\Vb,\bar{\eta},\bar{\mu})=\Mon(\V,\eta,\mu) }$}
   \put(650,180){\oval(100,100)[b]}
  \put(600,180){\line(0,1){200}}
  \put(700,180){\vector(0,1){200}}

%\V
  \put(1500,0){${(\V,\eta,\mu)}$}
   \put(1750,180){\oval(100,100)[b]}
  \put(1700,180){\line(0,1){200}}
  \put(1800,180){\vector(0,1){200}}
 \end{picture}
\end{center}
$\beta^*$ is the functor of composing $\beta$. It has a left adjoint $Lan_\beta$, the left Kan extension along $\beta$.
For $C\in Set^\B$ it is given by the coend formula
\[ Lan_\beta(C)(X) = \int^{n\in \F}  X^n\times C(n]     \]
The equivalences
\[ i_\F: Set^\F\lra \End,\hskip 1cm  i_\B: Set^\B\lra \An  \]
are defined by left Kan extensions that might be given by the following formulas
\[ i_\F(G)(X)=  \int^{n\in \F} X^n\times G(n],\hskip 1cm i_\B(C)(X)=\sum_{n\in \o} X^n\otimes_n C(n] \]
where $G\in Set^\F$ and $C\in Set^\B$.

Then the functor $i^a: \An \ra \End$ is just an inclusion and its right adjoint $(-)^a$ is given by the formula
\[ F^a(X) = \sum_{n\in\o} X^n\otimes_n F(n] \]
for $F\in \End$. $(-)^a$ is associating to functors and natural transformations their `analytic parts'.

Note that both $\An$ and $\End$ are strict monoidal categories with tensor given by composition, and $i^a$ is a strict monoidal functor. Thus its right adjoint $(-)^a$ has a unique lax monoidal structure making the adjunction $i^a\dashv (-)^a$ a monoidal adjunction. This in turn gives us a monoidal monad $(\V,\eta,\mu)$ on $\An$.

We have (cf. \cite{Z2}) a $2$-natural transformation $\cU$
\begin{center} \xext=1000 \yext=350
\begin{picture}(\xext,\yext)(\xoff,\yoff)
%1-cells
\putmorphism(0,250)(1,0)[\phantom{\MonCat}`\phantom{\bCat}`\Mon]{1000}{1}a
\putmorphism(0,150)(1,0)[\MonCat`\bCat`]{1000}{0}a
\putmorphism(0,50)(1,0)[\phantom{\MonCat}`\phantom{\bCat}`|-|]{1000}{1}b
%2-cell
\put(500,100){\makebox(300,100){$\Da\; \cU$}}
\end{picture}
\end{center}
where $\MonCat$ is the $2$-category of monoidal categories, lax monoidal functors, and monoidal transformations; $\Mon$ is the $2$-functor associating monoids to monoidal categories, $|-|$ is the forgetful functor forgetting the monoidal structure, and $\cU$ is a $2$-natural transformation whose component at a monoidal category $M$ is the forgetful functor from monoids in $M$ to the underlying category of $M$:  $\cU_M: \Mon(M) \ra |M|$.

Applying $\cU$ to the monoidal adjunction and $i^a\dashv (-)^a$ and monoidal monad $\V$ we get an adjunction between categories of monoids and a monad on $\Mon(\An)$. The unnamed arrow is $\Mon(i^a)$. But the monoids in $\End$ and $\An$ are monads and hence we get the left most adjunction $\cQ_f^m\dashv \cP_a^m$ that we were looking for together with the monad $(\Vb,\bar{\eta},\bar{\mu})$ on the category of analytic monads.

There are free monads on finitary functors (cf. \cite{Barr}) and free analytic monads on analytic functors (cf. \cite{Z1}). Therefore, the functors $\widehat{U}$ and $U$ have left adjoints $\widehat{F}$ and $F$, respectively. The adjunctions $F\dashv U$ and $\widehat{F}\dashv \widehat{U}$ induce monads $\M$ and $\widehat{\M}$, respectively. $\widehat{\M}$ is the finitary version of what is called `the monad for all monads' in \cite{Barr}. Putting this additional data to the above diagram and simplifying it at the same time we get a diagram
\begin{center} \xext=3000 \yext=1450
\begin{picture}(\xext,\yext)(\xoff,\yoff)
\setsqparms[0`0`0`0;1600`750]
 \putsquare(450,450)[\Mnd`\End`\AMnd`\An;```]

%horizontal
\putmorphism(450,1250)(1,0)[\phantom{\Mnd}`\phantom{\End}`\widehat{F}]{1600}{-1}a
\putmorphism(450,1150)(1,0)[\phantom{\Mnd}`\phantom{\End}`\widehat{U}]{1600}{1}b

\putmorphism(450,500)(1,0)[\phantom{\AMnd}`\phantom{\An}`F]{1600}{-1}a
\putmorphism(450,400)(1,0)[\phantom{\AMnd}`\phantom{\An}`U]{1600}{1}b

%vertical
 \putmorphism(500,1150)(0,-1)[\phantom{\Mnd}`\phantom{2\bCat_{\times}}`\cQ_f^m]{700}{1}r
 \putmorphism(400,1150)(0,-1)[\phantom{\Mnd}`\phantom{2\bCat_{\times}}`\cP_a^m]{700}{-1}l
 \putmorphism(2100,1150)(0,-1)[\phantom{\End}`\phantom{\An}`(-)^a]{700}{1}r
 \putmorphism(2000,1150)(0,-1)[\phantom{\End}`\phantom{\An}`i^a]{700}{-1}l

%\Vb
  \put(450,50){\oval(100,100)[b]}
  \put(400,50){\line(0,1){300}}
  \put(500,50){\vector(0,1){300}}
  \put(0,150){${(\Vb,\bar{\eta},\bar{\mu})}$}

%\V
  \put(2050,50){\oval(100,100)[b]}
  \put(2000,50){\line(0,1){300}}
  \put(2100,50){\vector(0,1){300}}
  \put(2150,150){${(\V,\eta,\mu)}$}

%\M
  \put(2550,450){\oval(100,100)[r]}
  \put(2550,500){\line(-1,0){400}}
  \put(2550,400){\vector(-1,0){400}}
  \put(2400,550){${(\M,\eta,\mu)}$}

%\M_hat
  \put(2550,1200){\oval(100,100)[r]}
  \put(2550,1250){\line(-1,0){400}}
  \put(2550,1150){\vector(-1,0){400}}
  \put(2400,1300){${(\widehat{\M},\widehat{\eta},\widehat{\mu})}$}

 \end{picture}
\end{center}
In the above diagram the square of the right adjoints commutes. Thus, the square of the left adjoint commutes as well. This shows in particular that the free monad on an analytic functor is analytic.

The monad $\Vb$ is a lift of a monad $\V$ to the category of $\M$-algebras $\AMnd$ and, by \cite{Beck}, we obtain

\begin{theorem}
The monad $\M$ for analytic monads distributes over the monad  $\V$ for finitary functors, i.e. we have a distributive law
\[ \lambda: \M\V \lra \V\M \]
The category of algebras of the composed monad $\V\M$ on $\AMnd$ is equivalent to the category $\Mnd$ of all finitary monads on $Set$. $\boxempty$
\end{theorem}

\vskip 2mm
{\em Remark.} We arrived at the above theorem with essentially no calculations at all. It has obvious positive aspects but it does not give an idea what the above distributive law is like. We shall present below explicit formulas how to calculate the values of some functors mentioned above and we shall also describe the coherence morphism $\varphi$ on the monoidal monad $\V$. This coherence morphism generates the distributive law $\lambda$. $\lambda$ is an analog of the distributive law of combing trees (cf. \cite{SZ1}). We think that it is possible that this $\lambda$ will eventually provide a formal explanation of the argument sketched in the proof of Theorem 14 in \cite{BD}.

\vskip 2mm

First we describe the adjunction $i^a\dashv (-)^a$. We shall drop the inclusion $i^a$ when possible. Let $A\in An$ and $G\in \End$ and $X$ be a set. The analytic functor $A$ is given by its coefficients. Its value at $X$ is
\[ A(X)=\sum_{n\in \o} X^n\otimes_n A_n\]
where $A_n$ is an $S_n$-set for $n\in \o$. The value of $G^a$ at $X$ is
\[ G^a(X)= \sum_{n\in \o} X^n\otimes_n G(n]\]
Thus
\[ \V(A)(X)=A^a(X) =\sum_{n,m\in \o} X^n\otimes_n(n]^m\otimes_m A_m  \]

The unit of the adjunction $i^a\dashv (-)^a$ at $X$
\[ (\eta_A)_X : A(X) \lra A^a(X) \]
is given by
\[ [\vec{x},a] \mapsto [\vec{x},1_n,a]   \]
where $\vec{x}:(n]\ra X$ and $a\in A_n$.

The counit of the adjunction at $X$
\[ (\varepsilon_G)_X : \sum_{n\in\o} X^n \otimes_n G(n]\lra G(X) \]
is given by
\[ [\vec{x},t] \mapsto G(\vec{x})(t)   \]
where $\vec{x}:(n]\ra X$ and $t\in G(n]$.

The multiplication in the monad $\V$
\[ (\mu_A)_X : \sum_{n,m,k\in \o} X^n\otimes_n(n]^m\otimes_m(m]^k\otimes_k A_k  \lra \sum_{n,k\in \o} X^n\otimes_n(n]^k\otimes_k A_k \]
is given by composition
\[ [\vec{x},g,f,a]\mapsto  [\vec{x},g\circ f,a]  \]
where $\vec{x}:(n]\ra X$, $g:(m]\ra (n]$, $f:(k]\ra (m]$, and $a\in A_k$.
This ends the definition of the monad $\V$.

Now we shall describe the monoidal structure on $\V$.

If $B$ is another analytic functor, the $n$-th coefficient of the composition $A\circ B$ is given by
\[ (A\circ B)_n =  \sum_{m,n_1,\ldots, n_m\in \o, \sum_{i=1}^m n_i =n} (S_n\times B_{n_1}\times\ldots \times B_{n_m}\times A_m)_{/\sim_n}  \]
where the equivalence relation $\sim_n$ is such that for $\sigma\in S_n$, $\sigma_i\in S_{n_i}$, $\tau\in S_m$, $b_i\in B_i$, for $i\in (m]$ and $a\in A_m$ we have

\[  \lk \sigma, \sigma_1\cdot b_1,\ldots, \sigma_m\cdot b_m, \tau\cdot a \rk \sim_n
\lk \sigma\circ (\lk\sigma_{1},\ldots, \sigma_{m}\rk\star\tau), b_{\tau(1)},\ldots, b_{\tau(m)}, a \rk  \]
where $\star$ is the composition in the operad of symmetries $\Sym$.

The $n$-th coefficient of  $\V(A)\circ \V(A)$ is given by
\[ (\V(A)\circ \V(A))_n =  \sum_{m,n_i,k_i\in \o, \sum_{i=1}^m n_i =n} (S_n\times (n_1]^{k_1} \otimes_{k_1} A_{k_1}\times\ldots \times (n_m]^{k_m} \otimes_{k_m} A_{k_m}\times A_m)_{/\sim_n}  \]
and the  $n$-th coefficient of  $\V(A^2)$ is given by
\[ (\V(A^2))_n =  \sum_{m,k,k_i\in \o, \sum_{i=1}^m k_i =k} (n]^k\times A_{k_1}\times\ldots \times A_{k_m} \times A_m \]
The coherence morphism $\varphi$ for $\V$ at the $n$-th coefficient of the functor $A$ is
\[ \varphi_n: (\V(A)\circ \V(A))_n  \lra (\V(A^2))_n  \]
is given by
\[  \lk \sigma, [\sigma_1,a_1],\ldots ,[\sigma_m,a_m], \tau, a \rk \mapsto
\lk \sigma\circ(\lk\sigma_{1},\ldots ,\sigma_{m}\rk\star\bigodot\tau), a_{\tau(1)}\ldots a_{\tau(m)}, a\rk   \]
Note that this map is well defined at the level of equivalence classes.

As the functor $(-)^a:\End\ra \An$ is monadic every finitary functor is a $\V$-algebra on an analytic functor. For $G$ in $\End$ the corresponding algebra map $\alpha_G$ at set $X$
\[ \alpha_G(X): \V((G)^a)(X)=\sum_{n,m\in\o} X^n\otimes_n (n]^m\otimes_m G(m) \lra \sum_{n\in\o} X^n\otimes_n G(n) = (G)^a(X) \]
is given by
\[ (\vec{x},f,t)\mapsto (\vec{x},G(f)(t)) \]
where $\vec{x}:(n]\ra X$, $f:(m]\ra(n]$, $t\in G(m)$.

%\newpage
\section{Equational theories vs Operads}

In this section we study the relations between equational theories  and operads, both symmetric and rigid. In particular, we shall describe the adjunction  $\cQ_f^{eo}\dashv \cP_a^{oe}$ and the properties of the embeddings
$\cE_s$ and $\cE_{ri}$.
\begin{center} \xext=1600 \yext=450
\begin{picture}(\xext,\yext)(\xoff,\yoff)
   \putmorphism(800,100)(1,0)[\SOp`\ET`]{800}{0}a
   \putmorphism(800,50)(1,0)[\phantom{\SOp}`\phantom{\ET}`\cQ_f^{eo}]{800}{-1}b
   \putmorphism(800,150)(1,0)[\phantom{\SOp}`\phantom{\ET}`\cP_a^{oe}]{800}{1}a
   \putmorphism(0,100)(1,0)[\RiOp`\phantom{\SOp}`\cP]{800}{1}a
   \put(700,360){$\cP_p^{oe}$}
   \put(0,200){\line(0,1){120}} \put(0,320){\line(1,0){1600}} \put(1600,320){\vector(0,-1){120}}
\end{picture}
\end{center}

\subsection*{The functor $\cP_a^{oe} : \SOp \ra \ET$}

We start by defining the functor $\cP_a^{oe}$. Let $\cO$ be a symmetric  operad. We define an equational theory $\cP_a^{oe}(\cO)=(L,A)$. As the set of $n$-ary function symbols we put $L_n=\cO_n$ for $n\in \o$. The set of axioms $A$ contains the following equations in context:
\begin{enumerate}
  \item $\iota(x_1)=x_1 :  \vec{x}^1$ where $\iota\in \cO_1$ is the unit of the operad $\cO$;
%  \item $\iota(f(x_1, \ldots , x_n)) = f(x_1, \ldots , x_n) =f(\iota(x_1), \ldots , \iota(x_n)) :  \vec{x}^n$ for all $f\in \cO_n$; $\iota\in \cO_1$ is the unit of the operad;
   \item $f(f_1(x_1, \ldots , x_{k_1}), \ldots ,f_m(x_{k_{m-1}+1}, \ldots , x_{k_m})) = (\lk f_1,\ldots, f_m\rk\ast f)(x_1, \ldots x_{k}) : \vec{x}^k$

   where $f\in \cO_m$, $f_i\in \cO_{k_i}$ for $i\in 1,\ldots, m$, $k=\sum_{i=1}^mk_i$;
     \item $ f(x_{\sigma(1)}, \ldots , x_{\sigma(n)})=(\sigma\cdot f)(x_1, \ldots , x_n) : \vec{x}^n$ for all $f\in \cO_n$ and $\sigma\in S_n$.
\end{enumerate}
Clearly, all equations are linear-regular and hence the theory $\cP_a^{oe}(\cO)$ is linear-regular.

Suppose that $h:\cO\ra \cO'$ is a morphism of symmetric operads. We define the interpretation $$\cP_a^{oe}(h): \cP_a^{oe}(\cO)\lra \cP_a^{oe}(\cO')$$ For $f\in \cO_n$ we put
\[ \cP_a^{oe}(h)(f) = (h(f)(x_1,\ldots, x_n):\vec{x}^n),\]
 for $n\in \o$.

\begin{proposition}\label{commuation_et_lt}
The following triangle
\begin{center} \xext=800 \yext=500
\begin{picture}(\xext,\yext)(\xoff,\yoff)
 \settriparms[1`-1`-1;400]
  \putVtriangle(0,0)[\ET`\LT`\SOp;\cL_e`\cP_a^{oe}`\cP_a]
\end{picture}
\end{center}
commutes up to a natural isomorphism.
\end{proposition}
{\em Proof.}  Let $\cO$ be a symmetric operad. We define a functor
\[ \psi_\cO : \cP_a(\cO) \lra \cL_e \cP_a^{oe}(\cO) \]
by
\[ [ \phi,!,f ] : n\ra 1 \mapsto [f(x_{\phi(1)},\ldots x_{\phi(m)}):\vec{x}^n] \]
where $\phi:(m]\ra (n]$, $f\in \cO_m$. The extension of this definition to morphisms with arbitrary codomains is obvious but it only complicates the notation.

$\psi_\cO$ is clearly bijective on objects.  Since every term in $\cP_a^{oe}(\cO)$ is provably equal to a simple term (=operation applied to variables), $\psi_\cO$ is full.

We shall show that $\psi_\cO$ is faithful. This is where combinatorics meets equational logic.
Suppose we have two morphisms $\lk \phi,!,g\rk$, $\lk \phi',!,g'\rk$ in $\cP_a(\cO)$
\begin{center} \xext=800 \yext=920
\begin{picture}(\xext,\yext)(\xoff,\yoff)
 \settriparms[1`1`0;400]
  \putAtriangle(0,400)[m`n`1;\phi`\lk !,g\rk`]
  \settriparms[0`-1`-1;400]
  \putVtriangle(0,0)[\phantom{n}`\phantom{1}`m';`\phi'`{\lk !,g'\rk}]
  \putmorphism(0,400)(1,0)[\phantom{n}`m`\bar{\phi}]{300}{-1}a
  \putmorphism(400,800)(-1,-4)[\phantom{n}`\phantom{m}`\sigma]{100}{1}r
  \putmorphism(370,120)(-1,4)[\phantom{n}`\phantom{m}`\sigma']{35}{1}r
\end{picture}
\end{center}
such that $\psi_\cO(\phi,!,g)=\psi_\cO(\phi',!,g')$. This means that the theory $\cP_a^{oe}\cO)$
proves
\[ g(x_{\phi(1)}, \ldots, x_{\phi(m)})= g'(x_{\phi'(1)}, \ldots, x_{\phi'(m')}) :\vec{x}^n \]

Since $\cP_a^{oe}(\cO)$ is linear-regular theory, $m=m'$ and there are permutations $\sigma,\sigma'\in S_m$
and a function $\bar{\phi} : (m]\ra (n]$ such that $\cP_a^{oe}(\cO)$ proves
\[ g(x_{\sigma(1)}, \ldots, x_{\sigma(m)})= g'(x_{\sigma'(1)}, \ldots, x_{\sigma'(m)}) :\vec{x}^m \]
and
\[ \phi=\bar{\phi}\circ \sigma,\;\;\; \phi'=\bar{\phi}\circ \sigma' \]
Thus  $\cP_a^{oe}(\cO)$ proves
\[ g(x_{1}, \ldots, x_{m})= g'(x_{\sigma^{-1}\sigma'(1)}, \ldots, x_{\sigma^{-1}\sigma'(m)}) :\vec{x}^m \]
and
\[ g(x_{1}, \ldots, x_{m})= (\sigma^{-1}\sigma')\cdot g'(x_1, \ldots, x_m) :\vec{x}^m \]
The last equality holds only if
\[ g= (\sigma^{-1}\sigma')\cdot g' \]
holds in $\cO$. But this together with $\phi'=\phi\circ \sigma^{-1}\circ \sigma'  $ means that
\[ (\phi,!,g)=(\phi',!,g') \]
in $\cP_a(\cO)$. Thus $\psi_\cO$ is faithful as well.
 $\boxempty$

Next we identify the image of the functor $\cP_a^{oe}$.

\begin{proposition}
The functor $\cP_a^{oe}$ is faithful, full on isomorphisms and its essential image is the category of linear-regular theories $\LrET$ i.e. it factorizes
as an equivalence of categories $\cE_o$ followed by $\cP_a^e$
\begin{center}
\xext=600 \yext=650
\begin{picture}(\xext,\yext)(\xoff,\yoff)
 \settriparms[-1`-1`-1;600]
 \putbtriangle(0,0)[\ET`\LrET`\SOp;\P_a^e`\cP_a^{oe}`\cE_o]
\end{picture}
\end{center}
\end{proposition}
{\em Proof.} As $\cL_e$ is an equivalence of categories, the fact that $\cP_a^{oe}$ is faithful and full on isomorphisms follows from Proposition \ref{commuation_et_lt} and the same properties of the functor $\cP_a$ stated in Proposition \ref{operad-emb}.

Let $I: \cP_a^{oe}(\cO)\lra \cP_a^{oe}(\cO')$ be a linear-regular interpretation. We shall define $h_{I} : \cO \lra \cO'$ such that $\cP_a^{oe}(h_I)=I$.
For $f\in \cO_n$, $I(f) : \vec{x}^n$ is a linear-regular term in $\cP_a^{oe}(\cO')$. As in $\cP_a^{oe}(\cO')$ every (linear-regular) term is provably equal to a simple (linear-regular) term (just one function symbol), we can assume that already $$I(f)=\bar{f}(x_{\sigma(1)},\ldots, x_{\sigma(n)}):\vec{x}^n$$ holds, where
$\bar{f}\in \cO'$. We put
\[ h_I(f)=\sigma\cdot \bar{f} \]
The verification that  $\cP_a^{oe}(h_I)=I$ is left for the reader.

Let $T=(L,A)$ be a linear-regular theory. We shall define a symmetric operad $\cO$ such that $T$ is isomorphic to $\cE_o(\cO)$. The set of $n$-ary operations $\cO_n$ is the set of linear-regular terms in context $\vec{x}^n$ modulo provable equations from the set of axioms $A$.  The group  $S_n$ acts of $\cO_n$ by permuting variables
\[ \sigma \cdot [t(x_1,\ldots,x_n):\vec{x}^n]= [t(x_{\sigma(1)},\ldots,x_{\sigma(n)}):\vec{x}^n] \]
The unit in $\cO_1$ is the term $[x_1:  \vec{x}^1]$. The composition in $\cO$ is defined by `disjoint substitution' i.e. before substituting terms we need to perform $\alpha$-conversion to make the result of the substitution a linear-regular term. For example substituting terms in contexts $[t_1(x_1,x_2):\vec{x}^2]$, $[t_2:\vec{x}^0]$ and   $[t_3(x_1,x_2,x_3):\vec{x}^3]$ into the term $[t(x_1,x_2,x_3):\vec{x}^3]$
we get
\[ [t(t_1(x_1,x_2),t_2,t_3(x_3,x_4,x_5)):\vec{x}^5]  \]
We hope that this explains the composition in $\cO$ better than a formal definition. It should be clear that $\cO$ is a symmetric operad.

There is an interpretation $I:T\ra \cP_a^{oe}(\cO)$ sending an operation $f\in L_n$ to the term in context
\[ [[f(x_1,\ldots,x_n):\vec{x}^n)] (\vec{x}^n):\vec{x}^n) ]  \]
Note that the term in context $[f(x_1,\ldots,x_n):\vec{x}^n)]$ is just a symbol of the theory $\cP_a^{oe}(\cO)$. There is also an interpretation  $I': \cP_a^{oe}(\cO)\ra T$ sending an operation $[f(x_1,\ldots,x_n):\vec{x}^n)]\in \cO_n$ to the same thing but considered this time a term in context $[f(x_1,\ldots,x_n):\vec{x}^n)]$ of the theory $T$. These two interpretations are mutually inverse. Thus $T$ is isomorphic to $\cP_a^{oe}(\cO)$ in $\ET$, as required.
$\boxempty$

\subsection*{The functor $\cP_p^{oe} : \RiOp \ra \ET$}

The functor $\cP^{oe}_{p}$ is defined as the composition of the functors $\cP^{oe}_a\circ \cP$. We have

\begin{proposition}
The functor $\cP_p^{oe}: \RiOp\ra \ET$ is faithful, full on isomorphisms and its essential image is the category of rigid theories $\RiET$.
\end{proposition}
{\em Proof.}
As $\cL_e:\ET\ra \LT$ is an equivalence of categories and $\cP:\RiOp\ra SOp$ is full and faithful, the fact that $\cP_p^{oe}$ is faithful and full on analytic morphisms (and hence also on isomorphisms) follows from Proposition \ref{commuation_et_lt} and the same properties of the functor $\cP^{ol}_p : \RiOp\ra \LT$ stated in Proposition \ref{im_rigid_in_lt}.

It remains to show that an equational theory is rigid iff it is of form $\cP_p^{oe}(\cO)$ for a symmetric operad $\cO$ whose actions on operations are free.
In the equational theory $\cP_p^{oe}(\cO)$ every term is equivalent to a function symbol $f\in \cO_n$ applied to some variables. Let us fix $n\in\o$ and   $f\in \cO_n$. Assume that for some $\sigma\in S_n$, the theory $\cP_p^{oe}(\cO)$ proves $$f(x_1,\ldots,x_n)=f(x_{\sigma(1)},\ldots,x_{\sigma(n)}):\vec{x}^n.$$
This means that in the Lawvere theory $\cP_p^{ol}(\cO)(n,1)$ we have equalities of analytic morphisms
\[ (id_n,!,f) = (\sigma,!,f) =(id_n,!,f)\circ \pi_\sigma . \]
But, by Proposition \ref{im_rigid_in_lt}, the actions of $S_n$ on analytic morphisms in $\cP_p^{ol}(\cO)(n,1)$ are free, i.e. $\sigma$ is the identity.
Since $f$ was arbitrary, $\cP_p^{oe}(\cO)$ is indeed a rigid equational theory. $\boxempty$

The following Corollary corrects a statement from \cite{CJ1} (cf. \cite{CJ2}) concerning characterization of equational theories corresponding to polynomial monads.

\begin{corollary}\label{eq_rigid-the_polymonad}
The equivalence of categories
\begin{center} \xext=1200 \yext=150
\begin{picture}(\xext,\yext)(\xoff,\yoff)
  \putmorphism(0,0)(1,0)[\ET`\Mnd`\cM_l\circ \cL_e]{1200}{1}a
\end{picture}
\end{center}
restricts to the equivalence between the category of rigid equational theories and the category of finitary polynomial monads on $Set$
\begin{center} \xext=1200 \yext=150
\begin{picture}(\xext,\yext)(\xoff,\yoff)
  \putmorphism(0,0)(1,0)[\RiET`\PMnd`\cM_e]{1200}{1}a
\end{picture}
\end{center}
\end{corollary}

%\newpage
\section{Comments and examples}

\begin{enumerate}
\item The equations expressing commutation of two operations are linear-regular. Therefore all operations in a theory $T$ commute iff they do in its analytic part $T^a$. However, the analytic part of an equational theory (or its monad on $Set$) is usually much bigger than the original equational theory. For example, if $\bT$ is a finitary monad on $Set$ then the value of its analytic part on one element set is the coproduct of the symmetrized free $\bT$ algebras on finitely many generators
      \[ \bT^a(1)= \sum_{n\in\o} 1^n\otimes_nT(n) = \sum_{n\in\o} T(n)_{/S_n}  \]
Thus it is not so surprising that theories arising in this way might be of interest only in special circumstances, preferably when the theory we start with is very small.

\item The categories $\SOp$ and $\LT$ are complete and cocomplete. The functor $\P_a:\SOp\ra \LT$ preserves all colimits as a left adjoint and it also preserves all connected limits. However, it does not preserve the terminal object. The terminal object is the value of $\cQ_f: \LT\ra \SOp$ on the terminal Lawvere theory. We describe it below.
\item Recall that ${\mathds{1}}$ denotes the terminal equational theory. It has one constant, say $e$, and can be axiomatized by a single axiom: $x_1=e:\vec{x}^1$. As a Lawvere theory it is the category that has exactly one morphism between any two objects. As $\cQ_f^e :\ET\ra \LrET$ is a right adjoint, it preserves the terminal object. Hence $\cQ_f^e(\mathds{1})$, the linear-regular part  ${\mathds{1}}$, is the terminal linear-regular theory. It is the theory of commutative monoids. It is best seen at the level of Lawvere theories. Both theories, $\cQ_f^e(\mathds{1})$ and the theory of commutative monoids, are linear-regular and, for any $n$, have exactly one analytic morphism
    \[a:n\ra 1\]
    In case of the theory of monoids it is given by
    \[ x_1, \ldots , x_n \mapsto  x_1\cdot \ldots \cdot x_n \]
  \item As we mentioned in Section 2, ${\mathds{1}}$ considered as a Lawvere theory has a proper subtheory,  in which $0\not\cong 1$. As equational theory it has no function symbols, and can be axiomatized by a single axiom: $x_1=x_2:\vec{x}^2$.  The analytic part of this theory is the theory of commutative semigroups.
  \item The embedding of the strongly regular theories into all equational theories has a right adjoint, $\cQ$, as well. The values $\cQ$ on the terminal equational theory  ${\mathds{1}}$ is the terminal strongly regular theory, i.e. the theory of monoids.
  \item As we saw above the Lawvere theory for monoids $\bT_{mon}$ is analytic. Thus it is an image under $\cL_o : \SOp\ra \ALT$ of a symmetric operad. It can be easily shown that any analytic morphism
      \[ a:n\ra 1 \]
      in $\bT_{mon}$ is of the form
      \[ x_1, \ldots , x_n \mapsto  x_{\sigma(1)}\cdot \ldots \cdot x_{\sigma(n)} \]
      where $\sigma\in S_n$, i.e. it is a multiplication of all variables in any order. Thus the operad $\bT^s_{mon}$ (see the proof of Proposition \ref{operad-emb-im} for notation $(-)^s$) is the operad of symmetries, $\Sym$, and hence the theory of monoids $\bT_{mon}$ is the image of the operad of symmetries under $\cL_o$.
  \item The Lawvere theory for monoids with anti-involution $\bT_{mai}$ is analytic, as well. Any analytic morphism \[ a:n\ra 1 \]  in $\bT_{mai}$ is of form
      \[ x_1, \ldots , x_n \mapsto  s^{\varepsilon_1}(x_{\sigma(1)})\cdot \ldots \cdot s^{\varepsilon_n}(x_{\sigma(n)}) \]
      where $\sigma\in S_n$, and $\varepsilon_i\in \{ 0,1\}$, for $i=1,\ldots, n$, and $s^0(x)=x$, $s^1(x)=s(x)$.
\end{enumerate}

\noindent
Instytut Matematyki,
Uniwersytet Warszawski,
ul. S.Banacha 2,
00-913 Warszawa, Poland\\
szawiel@mimuw.edu.pl,
zawado@mimuw.edu.pl

\end{document}